\documentclass[12pt]{article}
\usepackage{latexsym,amsfonts,amssymb}
\usepackage{csvsimple}
\usepackage{verbatim}
\usepackage{bm}
\usepackage{upgreek}
\usepackage[dvips]{color}
\usepackage{colordvi,multicol}
\usepackage{amsmath}
 \usepackage{graphicx}
 \usepackage[english]{babel}
\usepackage{caption}
\usepackage{tabu}
\usepackage{float}
\usepackage[pagewise, left, displaymath, mathlines]{lineno}
\usepackage{dsfont}
\usepackage[shortlabels]{enumitem}
\usepackage{algorithm}
\usepackage{algpseudocode}
\usepackage{algcompatible}
\usepackage{titlesec}
\usepackage{mathtools}
\usepackage{xcolor}
\titleformat{\paragraph}
{\normalfont\normalsize\bfseries}{\theparagraph}{1em}{}
\titlespacing*{\paragraph}
{0pt}{3.25ex plus 1ex minus .2ex}{1.5ex plus .2ex}

\newcommand{\ce}{C^\epsilon_t(1)}
\newcommand{\co}{C^0_t(1)}

\newcommand{\ces}{C^\epsilon_s(1)}
\newcommand{\cs}{C^0_s(1)}
\newcommand{\cm}{C^\epsilon_{s-}(1)}

\newcommand{\dint}{\displaystyle\int}
\def \cal{\mathcal}
\newcommand{\comm}[1]{}
\newtheorem{thm}{Theorem}[section]
\newtheorem{cor}[thm]{Corollary}
\newtheorem{lem}[thm]{Lemma}

\newtheorem{rem}[thm]{Remark}

\date{}

\begin{document}
\title{\bf Parameter estimation of stochastic SIR model driven by small L\'{e}vy noise with time-dependent periodic transmission}
\author{Terry Easlick  and Wei Sun\thanks{Corresponding author.}\\ \\
  {\small Centre de recherche Azrieli du CHU Sainte-Justine}\\
{\small D\'epartement de Math\'ematiques et de Statistique}\\
{\small Universit\'e de Montr\'eal}\\
{\small Montreal, Canada}\\ \\
 {\small Department of Mathematics and Statistics}\\
{\small Concordia University}\\
{\small Montreal, Canada}\\ \\
{\small  terry.easlick@umontreal.ca,\ \ \ \ wei.sun@concordia.ca}}

\maketitle

\begin{abstract}
\noindent We investigate the parameter estimation and prediction of two forms of
the stochastic SIR model driven by small L\'{e}vy noise with time-dependent periodic transmission. We present consistency and rate of convergence results for the least-squares estimators.
We include simulation studies using the method of projected gradient descent.
\end{abstract}



\section{Introduction}
The study of infectious diseases predates the Industrial Revolution with the work of John Graunt in the 17th century. The use of mathematics to study epidemiological phenomena dates back to the time of Daniel Bernoulli with his model of smallpox innoculation in the 18th century (cf. \cite{ephist}). Just as infectious diseases remain a persistent fact of life so does the rigorous study of them. Much of the modern work on epidemiological models begins with the classical Susceptible-Infected-Recovered (SIR) model introduced by Kermack and McKendrick \cite{ker} nearly a century ago and is defined as:
\begin{eqnarray*}
	\begin{cases}
		\frac{dX_t}{dt} =  - \beta X_tY_t, \\
		\frac{dY_t}{dt} = \left( \beta X_t - \gamma\right)Y_t,\\
		\frac{dZ_t}{dt} = \gamma Y_t,
	\end{cases}
	\end{eqnarray*}
	where $\beta$ and $\gamma$ are the transmission and recovery coefficients, respectively. Additionally, demographics may be introduced to include the birth coefficient $\Lambda$ and the mortality coefficient $\mu$ as:
	\begin{eqnarray*}
	\begin{cases}
		\frac{dX_t}{dt} = \Lambda - \mu X_t - \beta X_tY_t, \\
		\frac{dY_t}{dt} = \left[ \beta X_t - (\mu+\gamma )\right]Y_t,\\
		\frac{dZ_t}{dt} = \gamma Y_t - \mu Z_t.
	\end{cases}
	\end{eqnarray*}

Advances in mathematical methodology have allowed such deterministic models to be re-visited in various stochastic frameworks; hence, making for a more realistic proxy of reality (cf. e.g., \cite{T}, \cite{G}, \cite{JJS}, \cite{Bao}, \cite{Chen1}, \cite{Zhang0},  \cite{Zhou1}, \cite{Zhou2}, \cite{EL}, \cite{Liu}, \cite{Pri} and  \cite{EL2}).
More recently, the unified stochastic SIR (USSIR) model was introduced by the authors in \cite{eas}. The USSIR model is defined as
	\begin{eqnarray}\label{1March}
		\begin{cases}
		\begin{split}
		&dX_t =\ b_1(t,X_t,Y_t,Z_t)dt + \sum\limits_{j=1}^r \sigma_{1j}(t,X_t,Y_t,Z_t)dB_t^{(j)}\\
			&\indent\indent + \int_{\left\{\lvert u \rvert \leq 1\right\}}\hspace{-.3cm}H_1(t,X_{t-},Y_{t-},Z_{t-},u)\tilde{N}(dt,du)+ \int_{\left\{\lvert u \rvert > 1\right\}}\hspace{-.3cm}G_1(t,X_{t-},Y_{t-},Z_{t-},u)N(dt,du),
		\end{split}\\
		\begin{split}
		&dY_t =\ b_2(t,X_t,Y_t,Z_t)dt + \sum\limits_{j=1}^r\sigma_{2j}(t,X_t,Y_t,Z_t)dB_t^{(j)}\\
			&\indent\indent+ \int_{\left\{\lvert u \rvert \leq 1\right\}}\hspace{-.3cm}H_2(t,X_{t-},Y_{t-},Z_{t-},u)\tilde{N}(dt,du)+  \int_{\left\{\lvert u \rvert > 1\right\}}\hspace{-.3cm}G_2(t,X_{t-},Y_{t-},Z_{t-},u)N(dt,du),
		\end{split}\\
		\begin{split}
		&dZ_t =\ b_3(t,X_t,Y_t,Z_t)dt + \sum\limits_{j=1}^r\sigma_{3j}(t,X_t,Y_t,Z_t)dB_t^{(j)} \\
			&\indent\indent+\int_{\left\{\lvert u \rvert \leq 1\right\}}\hspace{-.3cm}H_3(t,X_{t-},Y_{t-},Z_{t-},u)\tilde{N}(dt,du) + \int_{\left\{\lvert u \rvert > 1\right\}}\hspace{-.3cm}G_3(t,X_{t-},Y_{t-},Z_{t-},u)N(dt,du).
		\end{split}\\
		\end{cases}
		\end{eqnarray}
	Hereafter, $\mathbb{R}_+$ denotes the set of all positive real numbers, $(B_t)_{t\ge 0}=(B^{(1)}_t,\dots,B^{(r)}_t)_{t \ge 0}$ is a standard
	$r$-dimensional Brownian motion, $N$ is a Poisson random measure on $\mathbb{R}_+\times (\mathbb{R}^l\backslash\{0\})$ with intensity
	measure $\mu$ satisfying $\int_{\mathbb{R}^l\backslash\{0\}}(1\wedge|u|^2)\mu(du)<\infty$ and $\tilde{N}(dt,du) = N(dt,du) - \mu(du)dt$,
	$(B_t)_{t \ge 0}$ and $N$ are independent, $b_i, \sigma_{ij}:[0,\infty) \times \mathbb{R}^3_+\mapsto\mathbb{R}$,
	$H_i,G_i:  [0,\infty) \times \mathbb{R}^3_+\times (\mathbb{R}^l\backslash\{0\}) \mapsto\mathbb{R}$, $i=1,2,3,\ j=1,2,\ldots,n$, are measurable functions.

The goal of this paper is to accomplish parameter estimation of a periodic transmission function present in a stochastic SIR model, e.g., the USSIR model (\ref{1March}). It is known that many phenomena are periodic in nature; thus, it is a natural question to ask if transmission of a disease is periodic. Understanding this periodicity greatly aids in ability to predict possible outcomes. 
A potential application is the study of periodicity during
the COVID-19 pandemic, as well as future pandemics. The stochastic model chosen has L\'{e}vy noise subjected
to a small coefficient $\varepsilon$. Inclusion of this coefficient ensures the noise is not unreasonably erratic. Namely, we expect the reality of disease spread and its measurement to be noisy but not outright chaotic.  We focus not only on estimating unknown parameters and predicting future
behaviour but also the asymptotics of our estimators.
Given the complexity of the model, we use optimization techniques to iteratively solve for approximations to the estimators.

In recent times, the study of parameter estimation in stochastic epidemiological
models has been gaining momentum, but there is still much work that needs to be completed. Notable works are available in the literature, see e.g.,  \cite{Feng}, \cite{Contact}, \cite{XLiu}, \cite{Mum1}, \cite{Pan}, \cite{Zhang}, \cite{Green},
\cite{Li}, \cite{Mum}, \cite{Jagan}, \cite{Wacker}, \cite{Alen}, \cite{Chen}, \cite{Sim}, \cite{Gir}, \cite{Paul},  \cite{Bodhi} and \cite{Kroger}. In our work, we complete our estimation without the assumption of having prior knowledge on the explicit form of the noise. Additionally, the periodicity we seek to estimate is not assumed to be tied to specific time periods such as days, weeks or months;
rather, whatever the observed time period may be is what determines the time scale in which the periodicity occurs.

In addition to papers on parameter estimation of epidemiological models, there is a growing number
of works which consider the more general setting of parameter estimation for stochastic differential equations (SDEs) with small L\'{e}vy
noises regardless of specific applications (cf. e.g., \cite{KA}, \cite{ms}, \cite{Uch}, \cite{Xu}, \cite{gl},
 \cite{Shi}, \cite{Sun}, \cite{DJ}, \cite{Long} and \cite{KO}).
The methodology in this paper is directly inspired by the results of  Long et al. in \cite{Sun} and \cite{Long}. More specifically,
we utilize the least-squares method in a manner similar to Long et al. in \cite{Sun} and \cite{Long}. This allowed us to estimate a true parameter $\theta_0$ of a discretely observed stochastic
process. We introduce a contrast function from which we are able to derive the least-squares estimators (LSEs) and obtain results on the consistency
and limiting distributions of the estimators.

We consider a multi-dimensional parameter $\theta$ and an associated stochastic process $(S^{\varepsilon}_t)_{t\ge 0}$ which satisfies the SDE:
\begin{eqnarray}\label{MarchE11}
	dS^{\varepsilon}_t&=&b(t,S^{\varepsilon}_t,\theta)dt+\varepsilon\Bigg\{\sigma(t,S^{\varepsilon}_t)dB_t\nonumber\\
	&&
	\quad\quad+{\int_{\{|u|\le1\}}}H(t,S^{\varepsilon}_{t-},u)\widetilde{N}
	(dt,du)+{\int_{\{|u|> 1\}}}G(t,S^{\varepsilon}_{t-},u)N
	(dt,du)\Bigg\},
	\end{eqnarray}
	where $0<\varepsilon<1$, $t\in[0,1]$, $S^{\varepsilon}_0=s\in \mathbb{R}^d$, $\theta\in \overline{\Theta}$, the closure of an open convex bounded subset $\Theta$ of $\mathbb{R}^p$,  $b(\cdot,\cdot,\cdot):[0,\infty)\times\mathbb{R}^d\times\Theta\to \mathbb{R}^d$,
	$\sigma(\cdot,\cdot):[0,\infty)\times \mathbb{R}^d \to \mathbb{R}^{d\times r}$,
	$H(\cdot,\cdot,\cdot)$, $G(\cdot,\cdot,\cdot):[0,\infty)\times \mathbb{R}^d \times  (\mathbb{R}^l\backslash\{0\}) \to \mathbb{R}^d$ are Borel measurable functions. Suppose
$(S^{\varepsilon}_t)_{t\ge 0}$ is observed at regularly
spaced time points $\{t_k = \frac{k}{n}$, $k = 1, 2,\dots, n\}$. Hereafter, we denote the transpose of a matrix $A$ by $A^\top$. Define the contrast function
  $$
  \Psi_{n,\varepsilon}(\theta)=\sum_{k=1}^nnP_k^\top (\theta)P_k(\theta),
  $$
  where
  $$
  P_k(\theta)= S^{\varepsilon}_{t_k} - S^{\varepsilon}_{t_{k-1}}
  -\frac{1}{n}b(t_{k-1},S^{\varepsilon}_{t_{k-1}}, \theta).
  $$
  Let $\hat{\theta}_{n,\varepsilon}$ be a minimum contrast estimator, i.e., a random variable satisfying
  $$
  \hat{\theta}_{n,\varepsilon} := \arg\min_{\theta\in\Theta}\Psi_{n,\varepsilon}(\theta).
  $$

As we will see for the USSIR model, finding a closed form for the LSE $\hat{\theta}_{n,\varepsilon}$ can be quite difficult, thus
we look for suitable approximations
$\hat{\theta}_{n,\varepsilon}^*$ of $\hat{\theta}_{n,\varepsilon}$. In the vernacular for optimization, we may refer to our contrast function $\Psi_{n,\varepsilon}(\theta)$ as an objective function. Given an objective function, gradient descent (GD) is a common method used to obtain minimization results (cf. \cite[Section 8.1]{beck} for more information). More explicitly, given a convex function  $f$ without
constraints such that the goal is to solve
$$ \arg\min_x f(x), $$
the method of GD aids in finding the minimizer by use of the update rule
$$
	x^{(k+1)} = x^{(k)} - \eta \nabla f(x^{(k)})
$$
for some initial value $x^{(0)}$ and learning rate $\eta > 0$.
With the presence of constraints, GD is still applicable in the form of projected gradient descent (PGD)  (cf. \cite[Section 10.2]{beck}). That is, given a convex function
$f$ and a constraint set $\mathcal{C}$, a similar update rule exists and is stated as
$$
 y^{(k+1)} = x^{(k)}  - \eta \nabla f(x^{(k)}) ;\ \  x^{(k+1)} = \arg\min_{x\in\mathcal{C}}|y^{(k+1)}-x|.
$$

The remainder of this paper is organized as follows. Section 2 contains the theoretical results on LSEs for time-dependent SDEs driven by small L\'{e}vy noises with applications to the USSIR model (\ref{1March}). Sections 3 and 4 cover simulation studies of stochastic SIR models for population numbers and population proportions, respectively. In section 5, we give the proofs of the results presented in section 2. Lastly, a closing discussion of our results is given with mention of future directions.

\section{Least-squares estimators for time-dependent SDEs driven by small L\'{e}vy noises}\setcounter{equation}{0}
\subsection{General time-dependent SDEs driven by small L\'{e}vy noises}

In this subsection, we investigate LSEs  for discretely observed stochastic
processes driven by small L\'evy noises. The results presented here generalize the results of Long et al. in \cite{Sun} and \cite{Long}. Our contributions are the generalization to
a time-dependent jump-diffusion model. This generalization includes singular and non-singular noise matrices related to the L\'{e}vy noise. The Brownian motion, small jump and large jump parts are permitted to have different coefficients rather than a shared coefficient as is the case in \cite{Long}. Moreover, we do not require any knowledge of the explicit form for the noise term beyond that is  assumed to be a L\'{e}vy noise -- which is itself very general. It is worth noting that a different contrast function is used compared to the contrast function in \cite{Long}. However, we also provide a result for the same contrast function given in \cite{Long}. As the primary focus of this paper is the parameter estimation of a stochastic SIR model, the proofs for the following results are left until Section 5.

Consider an underlying deterministic (ordinary) differential equation denoted as
$$
		dS^0_t = b(t,S_t^0,\theta_0)dt,\ t \in [0,1];\ \ \ \ S^0_0 = s,
$$
where $\theta_0$ is the true value of the drift parameter. Denote by $C_{\uparrow}^{1,k,l}([0,1]\times\mathbb{R}^d \times \Theta ; \mathbb{R}^d)$ the class of functions $f \in C^{1,k,l}([0,1]\times\mathbb{R}^d \times \Theta ; \mathbb{R}^d) $ which satisfy
$$
	\sup_{t \in [0,1]} \sup_{\theta \in \Theta} \lvert \partial^{\nu_3}_\theta\partial^{\nu_2}_x \partial^{\nu_1}_t f(t, x, \theta) \rvert
	\leq C(1+|x|)^{\lambda}
$$
for some constants $C$ and $\lambda$ where $\nu_1,\nu_2$, and $\nu_3$ are non-negative integer-valued multi-indices satisfying  $0 \leq \nu_1 \leq 1, \sum_{i=1}^d \nu_2^{(i)} \leq k  $ and $\sum_{i=1}^d \nu_3^{(i)} \leq l$.

We take the following assumptions, which are modifications of those given in \cite{Long}.

\noindent {\bf({A}1)}
For any $\varepsilon\in(0,1)$, the SDE (\ref{MarchE11}) admits a unique, strong solution $S^{\varepsilon}$ (cf. \cite[Theorem 2.2]{GS} for concrete sufficient conditions).

\noindent {\bf({A}2)}
There exist $K>0$ and $\eta,\xi\in {\cal B}_+(\mathbb{R}^l )$ such that for any $t\in[0,1]$, $x,y\in \mathbb{R}^d$, $u\in \mathbb{R}^l\backslash\{0\}$ and $\theta\in \overline{\Theta}$,
 \begin{eqnarray*}
 	&&|b(t,x,\theta) - b(t,y,\theta)| \le K|x-y|,\ \ \ \ \int_{\{|v|\le 1\}}\eta^2(v)\mu(dv)<\infty,\\
&& | b(t,x,\theta)|+|\sigma(t,x)|+\frac{1_{\{|u|\le 1\}}|H(t,x,u)|}{\eta(u)}+\frac{1_{\{|u|> 1\}}|G(t,x,u)|}{\xi(u)} \le K(1+|x|).
 \end{eqnarray*}

\noindent {\bf({A}3)}
$b(\cdot,\cdot,\cdot)\in C_{\uparrow}^{1,1,3}([0,1]\times\mathbb{R}^d \times \Theta ; \mathbb{R}^d)$.

\noindent {\bf({A}4)}  $\theta \neq \theta_0 \Leftrightarrow \exists t \in [0,1]$ such that  $b(t,S_t^0,\theta) \neq b(t,S_t^0,\theta_0)$.

\noindent {\bf({A}5)} $\sigma(\cdot,\cdot)$ is continuous on $[0,1]\times \mathbb{R}^d$ and $H(\cdot,\cdot,u)$, $G(\cdot,\cdot,u)$ are continuous on $[0,1]\times \mathbb{R}^d$ for any $u\in \mathbb{R}^l\backslash\{0\}$.

\noindent {\bf({A}6)}
  $\varepsilon=\varepsilon_n\rightarrow 0$ and $n\varepsilon \rightarrow \infty$ as $n \rightarrow\infty$.

 Due to the nature of the work herein, it  is not always possible to obtain a closed-form of the estimator $\hat{\theta}_{n,\varepsilon}$. Hence it makes sense to consider approximation, which is denoted by  $\hat{\theta}_{n,\varepsilon}^*$. We use the notation $o_p(1)$ for a sequence of
random vectors that converges to zero in probability and the notation $\xrightarrow{P_{\theta_0}}$ for convergence in probability under $P_{\theta_0}$. By virtue of \cite[Theorem 5.7]{van}, similar to \cite[Theorem 2.1]{Sun, Long}, we can prove the following result on the consistency of the LSEs.

  \begin{thm}\label{thm1} Let $\hat{\theta}^*_{n,\varepsilon}$ be any sequence of estimators with $\Psi_{n,\varepsilon}(\hat{\theta}^*_{n,\varepsilon})\le \Psi_{n,\varepsilon}(\theta_0)+o_p(1)$. Then, under conditions (A1)-(A4), we have
	\begin{equation*}
		\hat{\theta}^*_{n,\varepsilon} \xrightarrow{P_{\theta_0}} \theta_0 \text{ as } \varepsilon \to 0
		\text{ and } n \to \infty.
	\end{equation*}
	\end{thm}
	
	Define the matrix $I(\theta)=(I^{ij} (\theta) )_{1\le i,j\le p}$ by
	$$
	I^{ij}(\theta)=\int_0^1(\partial_{\theta_i}b)^\top(r,S^0_r,\theta)\partial_{\theta_j}b(r,S^0_r,\theta)dr.
	$$
Similar to \cite[Theorem 2.2]{Sun, Long}, we can prove the following result on the rate of convergence of the LSEs.	

	\begin{thm}\label{thm2} Assume that conditions (A1)-(A6) hold and $I(\theta_0)$ is positive definite. Then,
	\begin{eqnarray*}
	\varepsilon^{-1}(\hat{\theta}_{n,\varepsilon} -\theta_0)
	&\xrightarrow{P_{\theta_0}}&I^{-1}(\theta_0)\left(\int_0^1(\partial_{\theta_i}b)^\top(r,S^0_r,\theta)\Bigg\{\sigma(r,S^0_r)dB_r\right.\\
	&&\left.
	\ \ \ \ +{\int_{\{|u|\le1\}}}H(r,S^0_r,u)\widetilde{N}
	(dr,du)+{\int_{\{|u|> 1\}}}G(r,S^0_r,u)N
	(dr,du)\Bigg\}\right)^\top_{1\le i\le p}
	\end{eqnarray*}
	$\text{ as } \varepsilon \to 0
		\text{ and } n \to \infty$.
	\end{thm}
	
\begin{rem}\label{rem1} Consider the following SDE:
\begin{eqnarray*}
	dS^{\varepsilon}_t=b(t,S^{\varepsilon}_t,\theta)dt+\varepsilon\sigma(t,S^{\varepsilon}_t)\left\{dB_t
	+{\int_{\{|u|\le1\}}}u\widetilde{N}
	(dt,du)+{\int_{\{|u|> 1\}}}uN
	(dt,du)\right\}.
	\end{eqnarray*}
In the event the diffusion matrix $\sigma\sigma^\top$ is invertible, we may use the following contrast function from Long et al. \cite{Long}:
\begin{eqnarray}\label{March7a}
  \Psi_{n,\varepsilon}(\theta)=n\bigg(\sum_{k=1}^nP_k^\top(\theta)\Lambda_{k-1}^{-1}P_k(\theta)\bigg)1_{\{D>0\}},
\end{eqnarray}
  where
  $$
  P_k(\theta)= S^{\varepsilon}_{t_k} - S^{\varepsilon}_{t_{k-1}}
  -\frac{1}{n}b(t_{k-1},S^{\varepsilon}_{t_{k-1}}, \theta),\ \  \Lambda_{k-1} = [\sigma\sigma^\top](t_{k-1},S_{t_{k-1}}),\ \ D = \inf_{k=0,\ldots,n-1}\det \Lambda_k.
$$
Define the matrix $I(\theta)=(I^{ij} (\theta) )_{1\le i,j\le p}$ by
\begin{eqnarray}\label{March13a}
	I^{ij}(\theta)=\int_0^1(\partial_{\theta_i}b)^\top(r,S^0_r,\theta)[\sigma\sigma^\top]^{-1}(r,S^0_r)\partial_{\theta_j}b(r,S^0_r,\theta)dr.
\end{eqnarray}
We make the following additional assumption.

\noindent {\bf({A}7)}
There exists an open convex subset $\mathcal{U}\subset \mathbb{R}^d$ such that
	$S_t^0 \in \mathcal{U}$ for all $t \in [0,1]$, $\sigma$ is smooth on $[0,1]\times \mathcal{U}$, and $\sigma\sigma^\top$ is invertible on $[0,1]\times\mathcal{U}$.

 Following the arguments of \cite[Theorems 2.1 and 2.2]{Sun, Long}, we can prove the following result.

\begin{cor}\label{cor1} Assume that conditions ({A}1)-({A}7) hold and $I(\theta_0)$ defined by (\ref{March13a}) is positive definite. Then, the assertions of
	Theorems \ref{thm1} and \ref{thm2} hold for the LSEs derived from the contrast function (\ref{March7a}).
\end{cor}

\end{rem}

\subsection{Application to USSIR model with periodic transmission}

Let $d=3$. We use equation (\ref{MarchE11}) to write the USSIR model as
\begin{eqnarray*}
\begin{bmatrix}	
dX_t^\varepsilon \\
dY_t^\varepsilon \\
dZ_t^\varepsilon \\
\end{bmatrix}
&=&
\begin{bmatrix}
	b_1(t,X_t^\varepsilon,Y_t^\varepsilon,Z_t^\varepsilon, \theta) \\
	b_2(t,X_t^\varepsilon,Y_t^\varepsilon,Z_t^\varepsilon,\theta) \\
	b_3(t,X_t^\varepsilon,Y_t^\varepsilon,Z_t^\varepsilon,\theta)
\end{bmatrix} dt +\varepsilon\bigg\{\sigma(t,X_t^\varepsilon,Y_t^\varepsilon,Z_t^\varepsilon)dB_t \\ &&
\ \ \ \ \ \ \ \ +{\int_{\{|u|\le1\}}}H(t,X^{\varepsilon}_{t-},Y^{\varepsilon}_{t-},Z^{\varepsilon}_{t-},u)\widetilde{N}
(dt,du)\nonumber\\ &&\
 \ \ \ \ \ \ \ +{\int_{\{|u|> 1\}}}G(t,X^{\varepsilon}_{t-},Y^{\varepsilon}_{t-},Z^{\varepsilon}_{t-},u)N
(dt,du)\bigg\}.\nonumber	
\end{eqnarray*}
Here the drift function $b=(b_1,b_2,b_3)$ is given in more general terms but in the subsequent sections explicit instances will be
given. Suppose that all previously stated assumptions hold.

The USSIR model was constructed to be a very general framework which is capable of reflecting real-world phenomena; however, it is also useful in studying hypothetical or artificial  instances. Namely, we gave priority to the biological relevance of the model; however, we did not restrict the possibility of artificiality in the model. The aim is that our methodology and theoretical results may be tested in a general setting. For instance, it is sound reasoning that the noise component of the susceptible compartment may not involve effects from the recovered populations -- from a biological/epidemiological viewpoint. However, from a mathematical standpoint, one way to ensure there is no reduction of dimension (degeneration) of the model is to take the noise of any given compartment as being affected by all compartments (i.e., assigning the noise coefficient of any compartment as a function of all compartments: $X_t,Y_t,Z_t$). Hence, this example of artificiality ensures the generality while testing our methods; yet, maintains the potential for specificity in applications. For the interested reader, we refer them to our previous work on the USSIR model (cf. \cite{eas}).

The centre of our focus in this paper is the presence of a periodic transmission function. Note that
any well-behaved periodic function may be approximated using a Fourier series. We consider the drift function
to contain a periodic transmission function of the form
\begin{equation}\label{beta}
	\beta(t, \theta) = \alpha_0 + \sum_{k=1}^K \alpha_{1,k}\cos\left(\frac{2\pi kt}{\vartheta}\right)  + \alpha_{2,k}\sin\left(\frac{2\pi kt}{\vartheta}\right),
\end{equation}
where $\theta = (\vartheta,\alpha_0,(\alpha_{1,k},\alpha_{2,k})_{k=1}^K)$ $\in \mathbb{R}_+^{2K+2}$ for $\vartheta,\alpha_0,\alpha_{1,k},\alpha_{2,k}>0,1\leq k\leq K$. The contrast function is then denoted by  $\Psi_{n,\varepsilon}(\theta)$. Moreover, we assume that $\vartheta \in [0,1]$, a natural assumption
given the estimation happens from time $t=0$ to $t=1$. The contrast function $\Psi_{n,\varepsilon}(\theta)$ is not globally convex; however, for a fixed $\vartheta$, it becomes convex in the remaining parameters. Hence, we introduce the following algorithm:

\begin{algorithm}[H]\label{ssgd}
	\caption{\ \ Linear-Search Gradient-Descent (LS-GD)}
	{\bf For} {$i$ in $1$ to $M$,} \\
	\phantom{Step}Fix a test value of $\vartheta$: $\vartheta_i \in \left( \frac{i-1}{M}, \frac{i}{M}  \right)$;\\
	\phantom{Step}Set the initial value $(\alpha_0^{(0)},(\alpha_{1,k}^{(0)},\alpha_{2,k}^{(0)})_{k=1}^K)$;\\
	\phantom{Step}Run Gradient Descent on the function $\Psi_{n,\varepsilon}(\vartheta_i, \alpha_0,(\alpha_{1,k},\alpha_{2,k})_{k=1}^K)$ with update rule\\
	\phantom{Step}$(\alpha_0^{(l)},(\alpha_{1,k}^{(l+1)},\alpha_{2,k}^{(l+1)})_{k=1}^K) =
	(\alpha_0^{(l)},(\alpha_{1,k}^{(l)},\alpha_{2,k}^{(l)})_{k=1}^K) - \eta \nabla_{\{\alpha_0,(\alpha_{1,k},\alpha_{2,k})_{k=1}^K\}}\Psi_{n,\varepsilon}(\vartheta_i, \alpha_0^{(l)},(\alpha_{1,k}^{(l)},\alpha_{2,k}^{(l)})_{k=1}^K )$;\\
	\phantom{Step}Store $\Psi^*_i := \arg\hspace{-2em}\min\limits_{\{\alpha_0,(\alpha_{1,k},\alpha_{2,k})_{k=1}^K\}}\hspace{-1em} \Psi_{n,\varepsilon}(\vartheta_i, \alpha_0,(\alpha_{1,k},\alpha_{2,k})_{k=1}^K )$.\\
	{\bf End}\\
	Return $  \hat{\theta}^* := \min_{i \in \{1,\ldots, M\}}\Psi^*_i$.
\end{algorithm}

\noindent The returned value $ \hat{\theta}^* $ will be the approximation of the LSE $\hat{\theta}=(\hat{\vartheta},\hat{\alpha}_0,(\hat{\alpha}_{1,k},\hat{\alpha}_{2,k})_{k=1}^K)$ we seek.
\begin{rem}
	The above algorithm is modified for PGD without difficulty. That minor modification is the algorithm form used in the simulation studies. Namely, it is a box-constraint problem, e.g., $0 \leq \vartheta \leq 1$.
\end{rem}

For the rest of the paper, we take $K=1$. That is, we
consider the case when $\beta(t, \theta) = \alpha_0 + \alpha_{1}\cos(2\pi t/\vartheta)  + \alpha_{2}\sin(2\pi t/\vartheta)$.
Choosing $K>1$ does not contribute any more insights. Additionally, choosing $K>1$ would be more computationally expensive to evaluate in our simulation studies.


\section{Simulation study of  SIR model for population numbers}\setcounter{equation}{0}

For this study, the unknown parameters are represented by the vector $\theta = \left( \vartheta,\alpha_0,\alpha_1,\alpha_2 \right) \in \mathbb{R}^4_+$
and our objective function is denoted $\Psi_{n,\varepsilon}\left( \theta \right)$, where
$n$ is the number of observations we have and $\varepsilon \in (0,1)$.
Note that, for a fixed $\vartheta$,  $\Psi_{n,\varepsilon}\left( \theta \right)$ is convex with respect to $(\alpha_0,\alpha_1,\alpha_2)$.

For our simulation and numerical estimation results we make use of the Julia programming language. This language is well suited for
this application as it was created with the purpose of being a high-level scientific computing centric programming language.
We use two well maintained Julia packages:  {\it DifferentialEquations.jl} and {\it Optim.jl } (cf. \cite{optim} and \cite{DiffEq}). The package {\it DifferentialEquations.jl} allows us to generate synthetic data for our testing purposes -- we generate the data such that we have $100$ observations.
A jump-adapted Euler-Maruyama discretization is the method for obtaining the $100$ observations. We next construct the contrast function and approximate the LSE by use of the LS-GD algorithm.
The task of implementing the LS-GD algorithm is accomplished using {\it Optim.jl} using in-built methods for gradient descent and dealing with constraints.  Namely, we pass
the constraint set $\mathcal{C}$ as an argument to the suitable method in {\it Optim.jl}. For a more in-depth exploration, we refer the reader to the {\it DifferentialEquations.jl} and {\it Optim.jl } documentation.

The contrast function $\Psi_{n,\varepsilon}$ is non-linear in the unknown parameters, which is clear from the definition of $\beta(t, \theta)$. Moreover,
there is interplay or dependence of optimal estimators between $\alpha_i$ and $\vartheta$ for $i =1,2$; however, our methodology allows us to handle this dependence.


\begin{rem}
	An initial value $(\alpha_0^{(0)},\alpha_1^{(0)},\alpha_2^{(0)})  = (0.51,0.31,0.21)$ was chosen for no particular reason other than it was suitable for all implementations of the LS-GD algorithm. There is nothing significant about this value; hence, we have no reason to believe it had any impact on our findings.

For the test values of $\vartheta_i$, we chose them such that $\vartheta_i \sim \cal{U}(\frac{i-1}{M}, \frac{i}{M})$. Hereafter, $\cal{U}$ denotes the uniform distribution.
	As the value $M$ increases, it does not make a difference if we individually choose the value
	or allow it to be randomly assigned within each subinterval.
\end{rem}

\subsection{Model for population numbers}

In this section, we consider the following model for population
numbers:
\begin{eqnarray}\label{ex1}
\begin{bmatrix} dX_t \\ dY_t \\ dZ_t \end{bmatrix} &=& \begin{bmatrix} \Lambda - \mu X_t - \beta(t, \theta)X_tY_t \\ \beta(t, \theta)X_tY_t - (\mu + \gamma)Y_t \\ \gamma Y_t - \mu Z_t\end{bmatrix} dt\nonumber\\
&& +\,
	 \varepsilon\! \begin{bmatrix}\sigma X_{t-}Y_{t-}Z_{t-} & 0 & 0 \\ 0 &\sigma X_{t-}Y_{t-}Z_{t-} & 0 \\ 0 & 0 &\sigma X_{t-}Y_{t-}Z_{t-} \end{bmatrix} \begin{bmatrix} dL^{(1)}_t \\ dL^{(2)}_t \\ dL^{(3)}_t \end{bmatrix},
\end{eqnarray}
where
$$
\begin{bmatrix} L^{(1)}_t \\ L^{(2)}_t \\ L^{(3)}_t \end{bmatrix} = \begin{bmatrix} B^{(1)}_t + \dint_0^t\dint_{\{|u|>  0.1\}}uN^{(1)}(ds,du)\\ B^{(2)}_t + \dint_0^t\dint_{\{|u|>  0.1\}} uN^{(2)}(ds,du)\\ B^{(3)}_t + \dint_0^t\dint_{\{|u|>  0.1\}} uN^{(3)}(ds,du) \end{bmatrix},
$$
$$\beta(t, \theta) = \alpha_0 + \alpha_1\cos\left(\frac{2\pi t}{\vartheta}\right)  + \alpha_2\sin\left(\frac{2\pi t}{\vartheta}\right) ,\ \ \ \ \theta \in [0,1],\ \alpha_0,\alpha_1,\alpha_2 > 0,$$
$\Lambda,\mu,\gamma, \sigma>0$ are constants, and $N = N_1 + N_2 $ such that $N_1$, $N_2$ are independent Poisson random measures with respective intensity
measures $\mu_1$ and $\mu_2 $:
\begin{equation}\label{j1rev}
	\begin{cases}
		N_1 \sim \frac{2}{3}\lambda dt\mu_1(du) \text{ with } \mu_1 = \delta(\{-0.1,0.1,0.0\}), \\
		N_2 \sim \frac{1}{3}\lambda dt\mu_2(du) \text{ with } \mu_2 = \delta(\{0.0,-0.1,0.1\}),
	\end{cases}
\end{equation}
where $\lambda \sim \cal{U}(\{1,2,3,4\})$.

By  \cite[Theorem 3.1]{eas},  the SDE (\ref{ex1}) has a unique strong solution taking values in $\mathbb{R}^3_+$. We utilize the  contrast function (\ref{March7a}), which is similar to that given in \cite{Long} and takes the explicit form:
$$
	\Psi_{\varepsilon,n}(\theta) =100\varepsilon^{-2}
	\left(\sum_{k=1}^{100} P_{t_k}^\top(\theta)\Lambda^{-1}_{k-1}P_{t_k}(\theta) \right),
$$
where
$$
P_{t_k}(\theta) = \begin{bmatrix}X_{t_k}  - X_{t_{k-1}} - \frac{1}{100} \left[ \Lambda - \mu X_{t_{k-1}} - \beta({t_{k-1}}, \theta )X_{t_{k-1}}Y_{t_{k-1}} \right]\\ Y_{t_k}  - Y_{t_{k-1}} - \frac{1}{100}\left[ \beta({t_{k-1}}, \theta )X_{t_{k-1}}Y_{t_{k-1}} - (\mu + \gamma)Y_{t_{k-1}} \right] \\  Z_{t_k}  - Z_{t_{k-1}}  - \frac{1}{100}\left(\gamma Y_{t_{k-1}}  - \mu Z_{t_{k-1}}  \right) \end{bmatrix},
$$
and
$$
\Lambda_{k-1} = \begin{bmatrix}(\sigma X_{t_{k-1}} Y_{t_{k-1}} Z_{t_{k-1}} )^2 & 0 & 0 \\ 0 &(\sigma X_{t_{k-1}} Y_{t_{k-1}} Z_{t_{k-1}})^2  & 0 \\ 0 & 0 & (\sigma X_{t_{k-1}} Y_{t_{k-1}} Z_{t_{k-1}})^2  \end{bmatrix}.
$$

\subsection{Parameter estimation and prediction}
We generated the synthetic data $10000$ times and randomly chose $1000$ such data sets to perform estimation upon--recall each data set contains $100$ observations. Set $\Lambda = 0.018,\mu = 0.00042$, $\gamma = 0.07142$, $\sigma = 0.5$ and $(X_0,Y_0,Z_0) = (2.3,0.19,0.25)$ (taken to be in millions).
For the unknown parameters, the true values used in generating the synthetic data are randomly decided such that
\begin{equation}\label{true}
\vartheta \sim \cal{U}(0,1),\ \  \alpha_0 \sim \cal{U}(0.1,0.8),\ \  \alpha_1,\alpha_2 \sim \cal{U}\left(0,\frac{\alpha_0}{\sqrt{2}}\right).
\end{equation}
It is noteworthy that although the distributions of $\alpha_1,\alpha_2$ are dependent on
the distribution of $\alpha_0$, their specific values are not. This dependence is necessary given we assumed $\beta(t, \theta) \ge 0$.
Also note that for each generation of synthetic data, the jump parameter $\lambda$ is randomly chosen as described in (\ref{j1rev}).





In the figures below, a majority of points lie close to the diagonal, reinforcing the quality of the obtained estimators. Certainly some outliers are to be expected; but, it stands to reason that as observations increase so will the accuracy. Additionally, if we focus on the plot for $\hat{\vartheta}^*$, it is clear that there is a wider spread among values progressively as they approach the point $(1,1)$. This effect can be explained by the estimation occurring on the time scale $t =0$ to $1$, and thus for a period $\vartheta$ closer to $1$ than $0$, it will have a lower frequency thus making it necessary to have more measurements (observations) as compared to a period closer to $0$.

\begin{figure}[H]\label{fig1}
\begin{center}
\includegraphics[width=\textwidth]{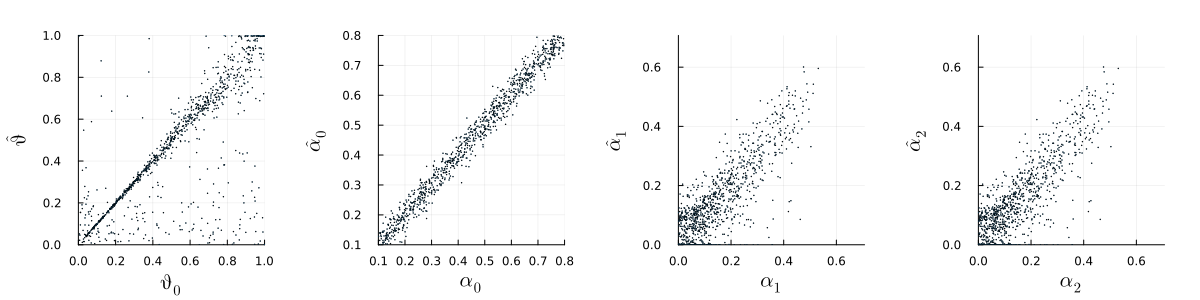}
\caption{Estimated values against true values  while  $\varepsilon = 0.3$. }
\end{center}
\end{figure}

\begin{figure}[H]\label{fig1}
	\begin{center}
	\includegraphics[width=\textwidth]{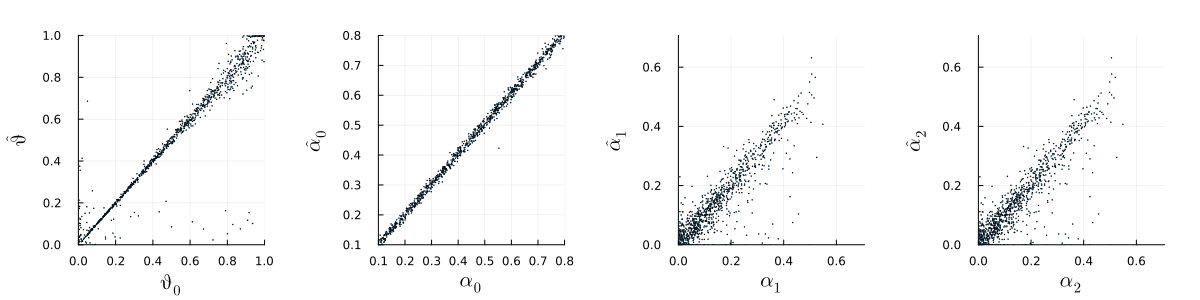}
	\caption{Estimated values against true values while  $\varepsilon = 0.1$. }
	\end{center}
	\end{figure}

\begin{figure}[H]\label{fig1}
		\begin{center}
		\includegraphics[width=\textwidth]{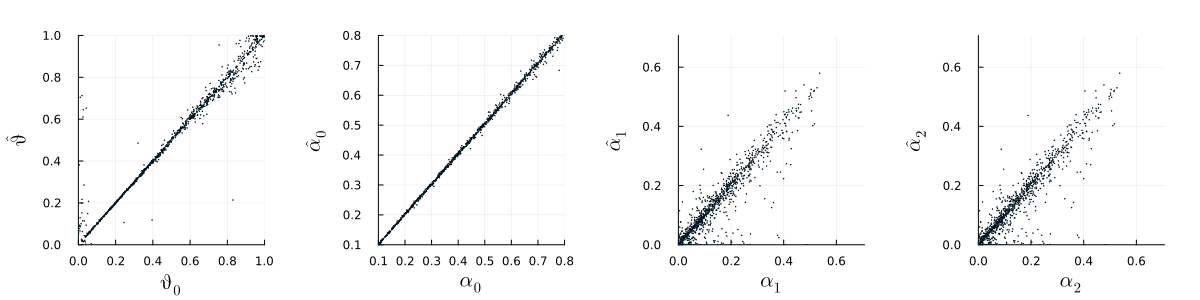}
		\caption{Estimated values against true values while  $\varepsilon = 0.01$. }
		\end{center}
\end{figure}	

\begin{figure}[H]\label{fig1}
\begin{center}
\includegraphics[width=\textwidth]{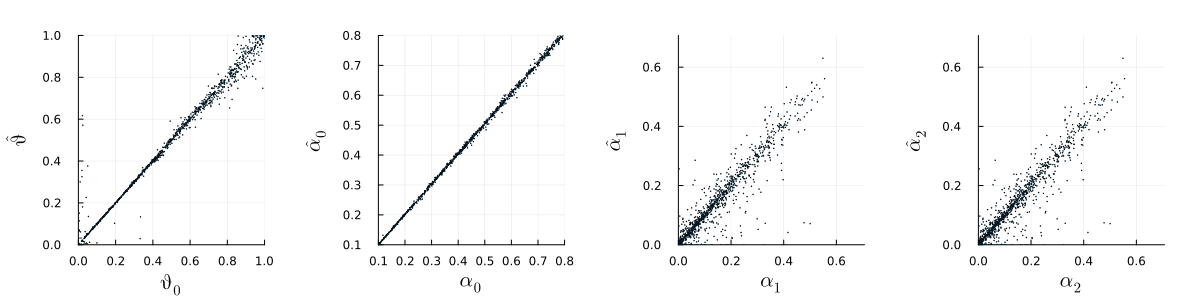}
\caption{Estimated values against true values while  $\varepsilon = 0.001$. }
\end{center}
\end{figure}

The previous results were applied to a simulation that goes beyond the estimation time in order to test the predictive ability. A true parameter $\theta_0$ is generated and used to create synthetic data. The initial conditions are randomly-generated to obtain the estimate $\hat{\theta}$ on the time interval $[0,1]$. Subsequently,  the true value $\theta_0$ and the estimate $\hat{\theta}$ are used to simulate the stochastic SIR model  (\ref{ex1}) on the time interval $[0,3]$ with newly generated initial conditions. For this purpose, we focused on the cases where $\varepsilon=0.3$ and $0.001$. The initial conditions were generated as $ (X_0, Y_0, Z_0) \sim (\cal{U}(1,4),\, \cal{U}(0.1,2),\, \cal{U}(0.1,1))$. The allowable ranges of values were chosen to reflect cases we found interesting. Namely, susceptible population is likely to be large so the disease has not completely overrun the population yet permitting the infected population to be unignorable in size.


The following table includes the randomly-generated true parameter values against the estimated values.

\begin{table}[H]
\centering
\begin{tabular}{|c | c | c | c|}
\hline
$(\vartheta, \alpha_0, \alpha_1, \alpha_2)$ & $(0.26836304, 0.15114833,0.0621514,0.096762)$ \\
\hline
$(\hat{\vartheta}, \hat{\alpha}_0, \hat{\alpha}_1, \hat{\alpha}_2)$, $\varepsilon = 0.3 $ & $(0.2759104,0.15106501,0.09268711, 0.0418802)$
\\
\hline
$(\hat{\vartheta}, \hat{\alpha}_0, \hat{\alpha}_1, \hat{\alpha}_2)$, $\varepsilon = 0.001 $ & $(0.26857489,0.151209307,0.06245031,0.0967899)$ \\	[.5ex]
\hline
\end{tabular}
\caption{True parameter values versus trained parameter values for model (\ref{ex1}).  }
\label{table:1}
\end{table}

The following plots contain the average trajectory from 100 simulations for the cases where $\varepsilon = 0.3$ and $0.1$ for the estimate $\hat{\theta}$. Whereas for the true parameter $\theta_0$ the path is the underlying deterministic path. In the plots, the path corresponding to the true parameter is denoted by $(X_t(\theta_0),Y_t(\theta_0),Z_t(\theta_0))$ and the estimated average trajectory by $(\overline{X}_t(\hat{\theta}_\varepsilon),\overline{Y}_t(\hat{\theta}_\varepsilon),\overline{Z}_t(\hat{\theta}_\varepsilon))$ where $\varepsilon = 0.3$ or $0.001$.

\begin{figure}[H]\label{fig1}
	\begin{center}
	\includegraphics[width=\textwidth]{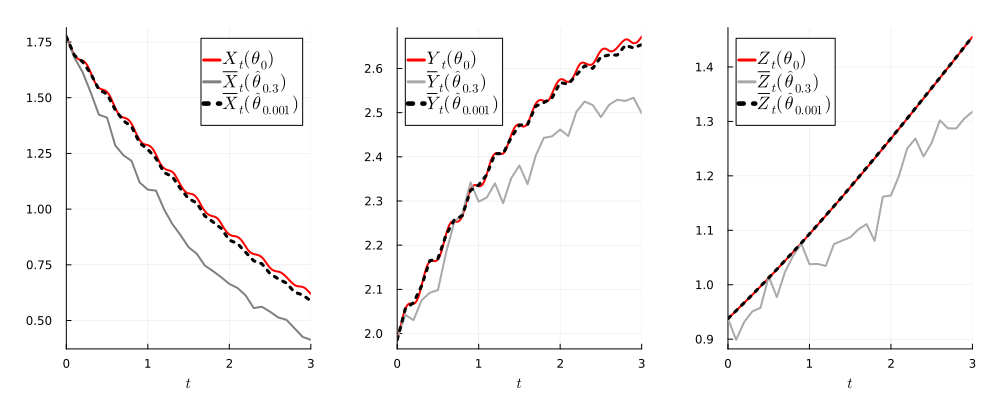}
	\caption{Simulation of model (\ref{ex1}) for  $(X_0,Y_0,Z_0)=(1.774, 1.985, 0.937)$. }
\end{center}
\end{figure}	

\begin{figure}[H]\label{fig1}
	\begin{center}
	\includegraphics[width=\textwidth]{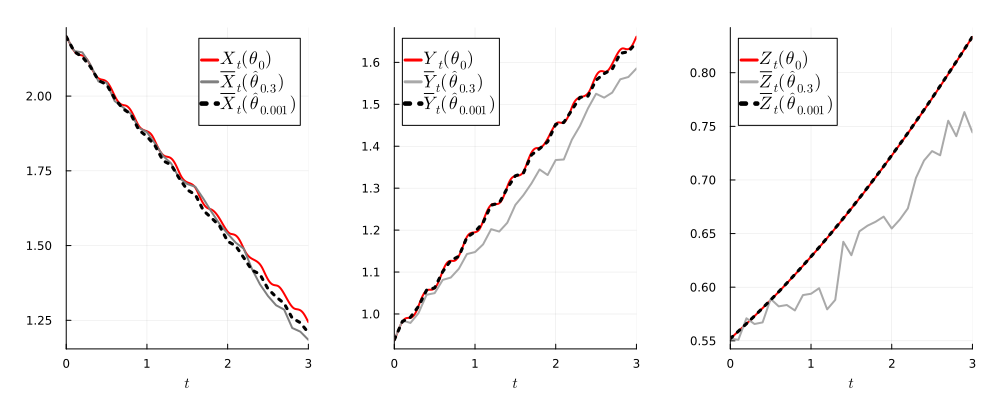}
	\caption{Simulation of model (\ref{ex1}) for  $(X_0,Y_0,Z_0)=(2.198, 0.938, 0.152)$.  }
	\end{center}
\end{figure}


\section{Simulation study of SIR model for population proportions}\setcounter{equation}{0}
We follow the same methodology of the previous section as applied to the population proportional model. The proportional model is a special case of the population numbers model (cf. \cite{eas}).
\subsection{Model for population proportions}
Consider the following model for population proportions:
\begin{equation}\label{ex2}
\begin{bmatrix} dX_t \\ dY_t \\ dZ_t \end{bmatrix} = \begin{bmatrix}  - \beta(t, \theta)X_tY_t \\ \beta(t, \theta)X_tY_t  - \gamma Y_t \\ \gamma Y_t \end{bmatrix} dt\, +\,
	 \varepsilon\! \begin{bmatrix}-\sigma X_{t-}Y_{t-}Z_{t-}   \\ 2\sigma X_{t-}Y_{t-}Z_{t-}   \\-\sigma X_{t-}Y_{t-}Z_{t-} \end{bmatrix} dL_t, \vspace{1cm}
\end{equation}
where
$$
L_t = B_t + \dint_0^t\dint_{\{|u|>  0.1\}} uN(ds,du),
$$
 and as before,
$$\beta(t, \theta) = \alpha_0 + \alpha_1\cos\left(\frac{2\pi t}{\vartheta}\right)  + \alpha_2\sin\left(\frac{2\pi t}{\vartheta}\right) ,\ \ \ \ \vartheta \in [0,1],\alpha_0,\alpha_1,\alpha_2 > 0,$$
$\gamma, \sigma>0$ are constants, and $N = N_1 + N_2$ as given in (\ref{j1rev}).

By  \cite[Theorem 2.1]{eas}, the SDE (\ref{ex2}) has a unique strong solution taking values
in the set $\Delta = \{(x,y,z) \in \mathbb{R}^3_+ : x+y+z =1\}$. It is easy to verify that all assumptions of Theorems 2.1 and 2.2 hold. We use the following contrast function:

$$
	\Psi_{\varepsilon,n}(\theta) =100\varepsilon^{-2} \left(\sum_{k=1}^{100} P_k^\top(\theta)P_k(\theta) \right),
$$
where
$$
P_k(\theta) = \begin{bmatrix}X_{t_k}  - X_{t_{k-1}} - \frac{1}{100} \left[- \beta({t_{k-1}}, \theta )X_{t_{k-1}}Y_{t_{k-1}} \right]\\ Y_{t_k}  - Y_{t_{k-1}} - \frac{1}{100}\left[ \beta({t_{k-1}}, \theta )X_{t_{k-1}}Y_{t_{k-1}} - \gamma Y_{t_{k-1}} \right] \\  Z_{t_k}  - Z_{t_{k-1}}  - \frac{1}{100}\gamma Y_{t_{k-1}}  \end{bmatrix}.
$$

\begin{rem}
It is noteworthy that the above contrast function differs from the one used for the population numbers model. Namely, this form does not require the noise matrix to be non-singular.
\end{rem}

\subsection{Parameter estimation and prediction}

To conform with the specifications of the proportional model given in (\ref{ex2}) we set $\gamma = 0.07142$, $\sigma = 0.5$ and $(X_0,Y_0,Z_0) = (0.82,0.07,0.11)$. Otherwise, the methodology is similar to that used in the previous section.
The following figures are presented in a similar fashion to the previous section with points lying on the diagonal signaling more accuracy in estimation.
\begin{figure}[H]\label{fig1}
	\begin{center}
	\includegraphics[width=\textwidth]{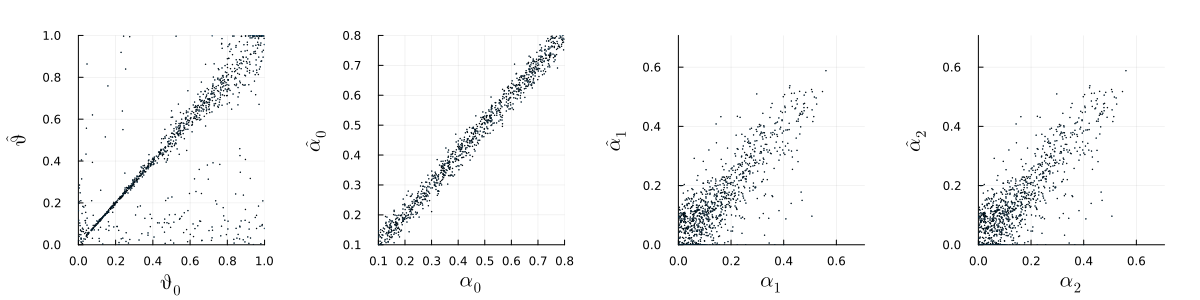}
	\caption{Estimated values against true values while $\varepsilon = 0.3$. }
	\end{center}
	\end{figure}
\begin{figure}[H]\label{fig1}
		\begin{center}
		\includegraphics[width=\textwidth]{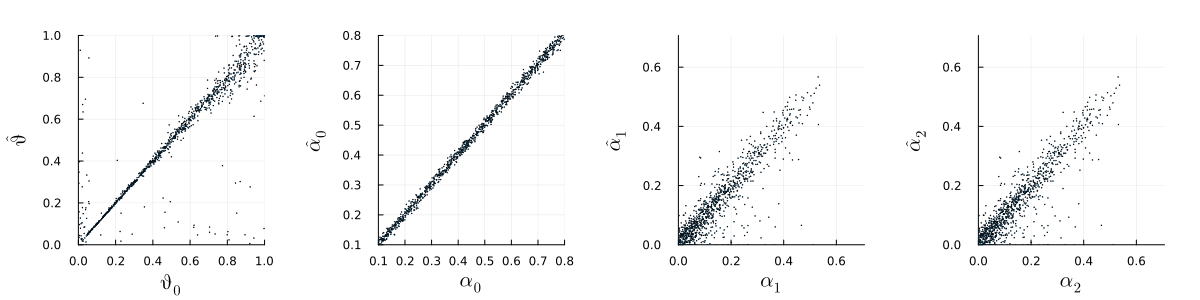}
		\caption{Estimated values against true values while $\varepsilon = 0.1$. }
		\end{center}
\end{figure}

\begin{figure}[H]\label{fig1}
	\begin{center}
	\includegraphics[width=\textwidth]{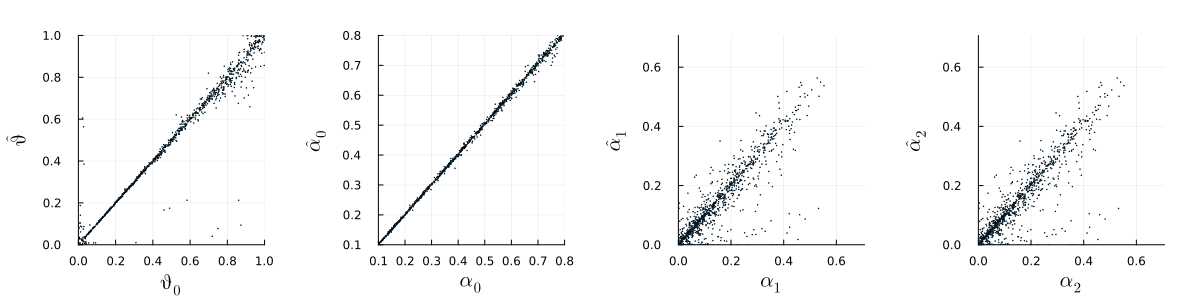}
	\caption{Estimated values against true values while $\varepsilon = 0.01$. }
	\end{center}
\end{figure}

\begin{figure}[H]\label{fig1}
	\begin{center}
	\includegraphics[width=\textwidth]{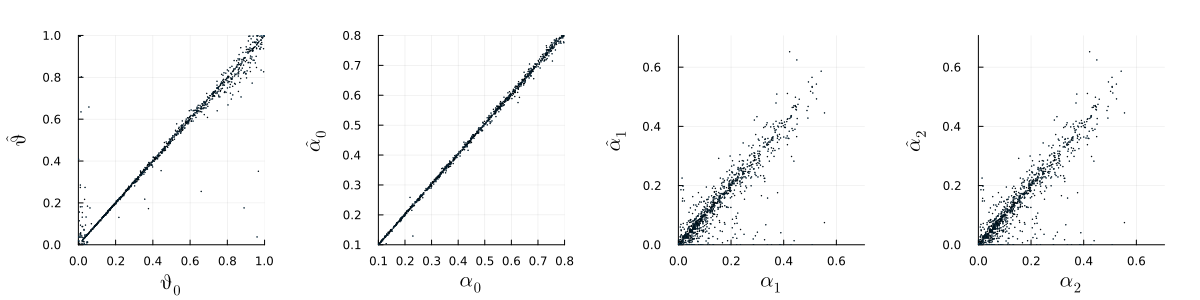}
	\caption{Estimated values against true values while $\varepsilon = 0.001$. }
	\end{center}
\end{figure}

The same methods as the previous section are used for this prediction study. We generate initial conditions such that each compartment in $(X_0,Y_0,Z_0)$ is sampled from $\cal{U}\left(0,1\right)$ with the constraint that $X_0 + Y_0 + Z_0 =1$ to conform to the proportional model definition.


The following table includes the true values against the trained (estimated) values for the proportional model.
\begin{table}[H]
	\centering
	\begin{tabular}{|c | c | c | c|}
	\hline
	$(\vartheta, \alpha_0, \alpha_1, \alpha_2)$ & $(0.26836304, 0.15114833,0.0621514,0.096762)$ \\
	\hline
	$(\hat{\vartheta}, \hat{\alpha}_0, \hat{\alpha}_1, \hat{\alpha}_2)$, $\varepsilon = 0.3 $ & $(0.27069755,0.1707693, 0.11016354, 0.082015210)$
	\\
	\hline
	$(\hat{\vartheta}, \hat{\alpha}_0, \hat{\alpha}_1, \hat{\alpha}_2)$, $\varepsilon = 0.001 $ & $(0.26743667816,0.151314344,0.058359562,0.0989769703)$ \\	[.5ex]
	\hline
	\end{tabular}
	\caption{True parameter values versus trained parameter values for proportional model.}
	\label{table:1}
	\end{table}

The subsequent plots are similar to those in the previous section in presentation. That is, they again display the average trajectories using the estimated parameters with the deterministic trajectory using the true parameter value.
\begin{figure}[H]\label{fig1}
	\begin{center}
	\includegraphics[width=\textwidth]{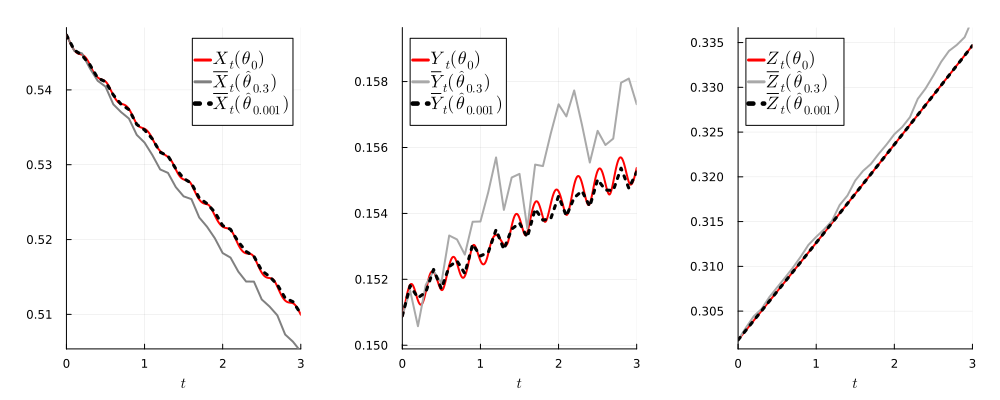}
	\caption{Simulation of model (\ref{ex2}) for such that $(X_0, Y_0, Z_0) =  (0.547, 0.152, 0.301) $. }
	\end{center}
\end{figure}	

\begin{figure}[H]\label{fig1}
	\begin{center}
	\includegraphics[width=\textwidth]{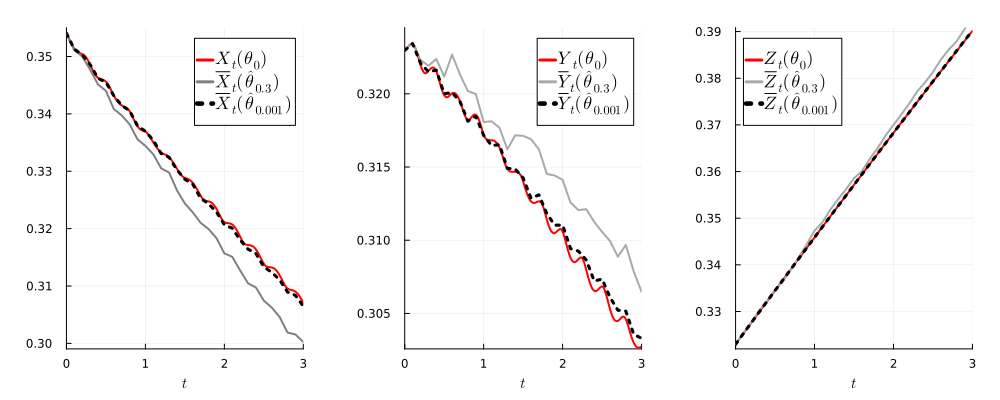}
	\caption{Simulation of model (\ref{ex2}) for such that $(X_0, Y_0, Z_0) =  (0.354, 0.324, 0.322) $.  }	\end{center}
\end{figure}


\section{Proofs of main results}\setcounter{equation}{0}

In this section, we give the proofs of Theorems \ref{thm1} and \ref{thm2}. Following the similar arguments of \cite{Long} with suitable modifications presented in this section,  Corollary \ref{cor1} can be proved. The details are omitted here.

Set $Y_t^{n,\varepsilon} = S^\varepsilon_{[nt]/n}$. To simplify notation, we use $C$ to denote a
generic constant whose value may vary from place to place, and use  $E$ and $P$ to denote $E_{\theta_0}$ and $P_{\theta_0}$, respectively.

\subsection{Proof of Theorem \ref{thm1}}

\begin{lem}\label{lem1}
Suppose {\bf (A1)} and {\bf(A2)} hold. Then, $\{Y_t^{n,\varepsilon}\}$ converges uniformly on compacts in probability to the deterministic process
$\{S_t^0\}$ as $\varepsilon\to0$ and $n\to\infty$.
\end{lem}
{\bf Proof.}\ \ We follow the argument of \cite[Lemma 4.1]{Long}. Let
$$
	J_t := N([0,t]\times \{|u| > 1\}),
$$
which  is a Poisson process with intensity $\lambda = \mu(\{|u| > 1\})$. Denote by
$\tau_1<\tau_2<\cdots<\tau_n<\cdots$ the jump times of $\{J_t\}$. We have $\lim_{n\rightarrow\infty}\tau_n=\infty$ a.s.. We use the interlacing technique (cf. \cite{Applebaum})
to construct the solution $(S^{\varepsilon}_t)_{t\ge 0}$ of (\ref{MarchE11}). Set $\tau_0=0$. Let  $\{Z_t^\varepsilon(i):t\ge 0\}$ be
the unique strong solution to the SDE:
\begin{eqnarray*}
&&Z^{\varepsilon}_t(i)= Z_0^\varepsilon(i) + \dint_0^t b(s,Z_s^\varepsilon (i),\theta)ds+\varepsilon\dint_0^t \left\{\sigma(s,Z_s^\varepsilon (i))dB_s(i)
+{\int_{\{|u|\le1\}}}H(s,Z_{s-}^\varepsilon (i),u)\widetilde{N}_i
(ds,du)\right\},\\
&&Z_0^\varepsilon(i)=S^{\varepsilon}_{\tau_{i-1}},
\end{eqnarray*}
where $B_t(i) = B_{\tau_{i-1}+t}-B_{\tau_{i-1}}$ and $N_i([0,t]\times A) = N([\tau_{i-1},\tau_{i-1}+t]\times A)$ for any $A \in\mathcal{B}(\mathbb{R}^l \backslash \{0\})$. Further, let
$\{Z_t^0 : t \ge 0\}$ be the unique strong solution to the underlying deterministic differential equation:
$$ Z_t^0(i) = Z_0^0(i) + \dint_0^t b(s,Z_s^0 (i),\theta)ds,\ \ \  Z_0^0(i) = S_{\tau_{i-1}}^0.$$
By {\bf (A1)}, we find that
$$
	Z_t^\varepsilon(i)=\left\{
	\begin{array}{ll}
		S^\varepsilon_{\tau_{i-1}+t},\ \ &0\leq t < \tau_i - \tau_{i-1}, \\
		S^\varepsilon_{\tau_{i-}},\ \ &t = \tau_i - \tau_{i-1},
	\end{array}\right.
$$
and $S^\varepsilon_{\tau_{i}} = S^\varepsilon_{\tau_{i-}} + \varepsilon G(\tau_{i-},S^\varepsilon_{\tau_{i-}},\xi_i)$, where $\{\xi_i : i \in \mathbb{N}\}$ are i.i.d. $\mathbb{R}^l$-valued random variables with common probability distribution $\frac{\mu(\cdot \cap \{|u|>1\})}{\lambda}$. First, we take $i=1$. Then, $Z_0^\varepsilon(1)=x\in \mathbb{R}^d$. Let $\varepsilon>0$ be fixed,
$$
	\overline{\tau}^\varepsilon_M := \inf\{t: |Z_t^\varepsilon(1)|\vee|Z_{t-}^\varepsilon(1)|>M\},
$$
and
$$
	f\in C_b^2(\mathbb{R}^d)\text{ such that } f(x) = |x|^2\text{ if } |x|\leq M.
$$
By  It\^o's formula (cf. \cite{Pro}),  we find that
$$
	f(Z_{t\wedge \overline{\tau}_M^\varepsilon}^\varepsilon(1)) - f(x) -\dint_0^{t\wedge \overline{\tau}^\varepsilon_M
	} A_sf(Z_{s}^\varepsilon(1))ds
$$
is a martingale, where
\begin{eqnarray*}
A_sf(x) &=& \sum_{k=1}^d b_i(s,x,\theta)\frac{\partial f}{\partial x_k}(x) +\frac{1}{2}\varepsilon^2\sum_{k,j=1}^d\sum_{l=1}^r
	\sigma_l^k(s,x)\sigma_l^j(s,x)\frac{\partial^2 f}{\partial x_k \partial x_j}(x) \\
	&&+\int_{\{|u| \le 1\}}\left [f(x +\varepsilon H(s,x,u)) - f(x) - \varepsilon\sum_{k=1}^d H_k(s,x,u) \frac{\partial f}{\partial x_k}(x)\right]\mu(du).
\end{eqnarray*}

Now we consider
\begin{eqnarray}\label{zeqn}
	Z_t^\varepsilon(1)-Z_t^0(1) &=& \dint_0^t [b(s,Z_s^\varepsilon(1),\theta)-b(s,Z_s^0(1),\theta)]ds + \varepsilon\dint_0^t
	\sigma(s,Z_{s}^\varepsilon(1))dB_s(1) \nonumber\\
	&&+\varepsilon \dint_0^t\dint_{\{|u| \le 1\}} H(s,Z_{s-}^\varepsilon(1),u)\tilde{N}_1(ds,du).
\end{eqnarray}
By Doob's inequality,  {\bf(A2)} and  (\ref{zeqn}), we get
\begin{eqnarray*}	
	&&{E}\left[\sup_{0 \leq s \leq t}\lvert Z_s^\varepsilon(1) - Z_s^0(1) \rvert^2\right]\\
&\le&4{E}[\lvert Z_t^\varepsilon(1) - Z_t^0(1) \rvert^2] \nonumber\\
& \leq&12
	{E}\Bigg\{ \bigg \lvert \dint_0^t b(s,\ces,\theta) - b(s,\cs,\theta)ds\bigg \rvert^2\nonumber\\&&\ \ \ \ \ \ \ \ 	+
	\varepsilon^2 \bigg \lvert\dint_0^t \sigma(s,\ces)dB_s(1)\bigg \rvert^2 +\varepsilon^2 \bigg \lvert\dint_0^t\dint_{\{|u| \leq 1\}} H(s,\cm,u)\tilde{N_1}(ds,du)\bigg \rvert^2 \Bigg\}\nonumber\\
		&\leq& 12K^2\Bigg\{t{E}\bigg [ \dint_0^t \lvert\ces - \cs \rvert^2ds\bigg] +
	\varepsilon^2 {E}\bigg [ \dint_0^t(1+ |\ces|)^2ds\bigg ]\nonumber\\
&&\ \ \ \ \ \ \ \ 	+\varepsilon^2 {E}\bigg [\dint_0^t\dint_{\{|u| \leq 1\}} \eta^2(u)(1+ |\cm|)^2\mu(du)ds\bigg]\Bigg\} \nonumber\\
&\leq& 12K^2\Bigg\{\left(t+2\varepsilon^2+2\varepsilon^2\int_{\{|u|\le 1\}}\eta^2(u)\mu(du)\right){E}\bigg [ \dint_0^t \lvert\ces - \cs \rvert^2ds\bigg]\\
&&\ \ \ \ \ \ \ \ \	+
	2t\varepsilon^2 (1+|x|)^2\left[1+\int_{\{|u|\le 1\}}\eta^2(u)\mu(du)\right]\Bigg\}.
\end{eqnarray*}
Then,  by Gronwall's inequality, we obtain that
\begin{eqnarray}\label{May28a}
	{E}\left[\sup_{0 \leq t \leq 1}\lvert \ce - \co \rvert^2\right] \to 0 \text{ as } \varepsilon \to 0.
\end{eqnarray}

For any small $\delta >0$, there is a $T$ sufficiently large enough such that ${P}(\tau_1 > T)<\delta$. Recall that the processes
$(S_t^\varepsilon)_{t\ge 0}$ and $(S_t^0)_{t\ge 0}$ are c\`{a}dl\`{a}g and continuous, respectively. Then, by (\ref{May28a}), we get
\begin{eqnarray*}	
	{P}\left(\sup_{0\leq t\leq \tau_1}|S_{t-}^\varepsilon - S_t^0|> \gamma\right) &\leq& {P}\left(\sup_{0\leq t\leq \tau_1}|S_{t-}^\varepsilon -S_t^0|>\gamma; \tau_1\leq T\right) +
	{P}(\tau_1 > T)\\
	&\leq&{P}\left(\sup_{0\leq t\leq T}|Z_{t}^\varepsilon(1) - Z_t^0(1)|>\gamma\right) + \delta\\
&\rightarrow&0\text{ as }\ \ \ \  \varepsilon,\,\delta \to 0.
\end{eqnarray*}
Additionally, note that $S^\varepsilon_{\tau_{1}} = S^\varepsilon_{\tau_{1-}} + \varepsilon G(\tau_{1-},S^\varepsilon_{\tau_{1-}},\xi_1)$. Since
$S_{\tau_1-}^\varepsilon \xrightarrow{P} S_{\tau_1}^0$, we get
$$
	|S_{\tau_1}^\varepsilon-S_{\tau_1}^0| \leq |S_{\tau_1-}^\varepsilon-S_{\tau_1}^0| + \varepsilon K(1+|S^\varepsilon_{\tau_1-}|)|\xi_1| \xrightarrow{P} 0.
$$
Thus,
$$
	{P}\left(\sup_{0\leq t\leq \tau_1}|S_{t}^\varepsilon -S_t^0|> \gamma\right) \to 0 \text{ as } \varepsilon\to 0.
$$

Next, we consider $\{Z_t^\varepsilon(2) : t\ge 0\}$, where $Z_0^\varepsilon(2) = S_{\tau_1}^\varepsilon$. We have
\begin{eqnarray*}	
	Z_t^\varepsilon(2)-Z_t^0(2) &=&S_{\tau_1}^\varepsilon - S_{\tau_1}^0 +
	\dint_0^t [b(s,Z_s^\varepsilon(2),\theta)-b(s,Z_s^0(2),\theta)]ds  \\
&& +\varepsilon\dint_0^t
	\sigma(s,Z_{s}^\varepsilon(2))dB_s(2) +
	\varepsilon \dint_0^t\dint_{\{|u| \le 1\}} H(s,Z_{s-}^\varepsilon(2),u)\tilde{N}_2(ds,du).
\end{eqnarray*}
Using the arguments above here yields
$$
	\lim_{\varepsilon\to0}{E}\bigg[\sup_{0\leq t\leq T}
		|Z_t^\varepsilon(2)-Z_t^0(2)|{1}_{\{|S_{\tau_1}^\varepsilon|+|S_{\tau_1}^0| \leq M\}} \bigg] = 0.
$$
For any small $\delta >0$, there exists $M > 0$ such that ${P}(|S_{\tau_1}^0|>\frac{M}{4})< \delta$.
Then, we have that
\begin{eqnarray*}	
&&\limsup_{\varepsilon \to 0} {P}\bigg(\sup_{0\leq t \leq T\wedge \tau_2} |Z_t^\varepsilon(2)-Z_t^0(2)|>\gamma \bigg)\\
&\leq& \limsup_{\varepsilon \to 0} {P}\bigg(\sup_{0\leq t \leq T} |Z_t^\varepsilon(2)-Z_t^0(2)|>\gamma; |S_{\tau_1}^\varepsilon|+|S_{\tau_1}^0| \leq M \bigg)\\
&&
	+ \limsup_{\varepsilon \to 0} {P}\bigg(|S_{\tau_1}^\varepsilon-S_{\tau_1}^0| > \frac{M}{2} \bigg)
	+ {P}\left(|S_{\tau_1}^0|>\frac{M}{4}\right)\\
&<&\delta.
\end{eqnarray*}
Note that $S^\varepsilon_{\tau_{2}} = S^\varepsilon_{\tau_{2-}} + \varepsilon G(\tau_{2-},S^\varepsilon_{\tau_{2-}},\xi_2)$. Thus,  we have
$$
	{P}\bigg( \sup_{\tau_1 \leq t \leq \tau_2} |S_t^\varepsilon - S_t^0| > \gamma \bigg) \to 0 \text { as } \varepsilon \to 0.
$$

By induction on $i\ge 1$, we get
$$
	{P}\bigg( \sup_{\tau_{i-1} \leq t \leq \tau_{i}} |S_t^\varepsilon - S_t^0| > \gamma \bigg) \to 0 \text { as } \varepsilon \to 0.
$$
As $i\to\infty$, it happens to be $\tau_i\to\infty$ a.s.,  thus for any given $T>0,\delta>0$ there exists $i_0 \in \mathbb{N}$ such that when
$i \ge i_0$ implies ${P}(\tau_i < T) < \frac{\delta}{2}$. Moreover,  for every $ i>0,\delta>0$, there exists $\varepsilon_0$ such that
for all $\varepsilon < \varepsilon_0$
$$
	{P}\bigg( \sup_{\tau_{i-1}\leq t \leq \tau_i} |S_t^\varepsilon - S_t^0|>\gamma\bigg)<\frac{\delta}{2^{i+1}i_0 }.
$$
Then, for $\varepsilon < \varepsilon_0$,
$$
	{P}\bigg(\sup_{0\leq t \leq 1} |S_t^\varepsilon - S_t^0| > \gamma \bigg) \leq
	\sum_{k=1}^i {P}\bigg(\sup_{\tau_{k-1}\leq t \leq \tau_k} |S_t^\varepsilon - S_t^0| > \gamma\bigg) < \frac{\delta}{2i_0}.
$$
Thus,
\begin{eqnarray*}	
	&&{P}\bigg(\sup_{0\leq t \leq 1} |S_t^\varepsilon - S_t^0| > \gamma\bigg)\\
&=& {P}\bigg(\sup_{0\leq t \leq 1} |S_t^\varepsilon - S_t^0| > \gamma; \tau_{i_0} < 1\bigg)+{P}\bigg(\sup_{0\leq t \leq 1} |S_t^\varepsilon - S_t^0| > \gamma; \tau_{i_0} \ge 1\bigg)\\
& \leq& {P}\bigg( \tau_{i_0} < 1\bigg) + \sum_{k=1}^{i_0} {P}\bigg(\sup_{0\leq t \leq 1}|S_t^\varepsilon-S_t^0|;\tau_{k-1} \leq t < \tau_k \bigg)\\
& <& \delta.
\end{eqnarray*}
Since $S_t^\varepsilon$ is a semimartingale, by \cite[Theorem II.11]{Pro}, we find that for any fixed $\varepsilon>0$,
$$
	\sup_{0\leq t \leq 1}|Y_{t}^{n,\varepsilon} - S_t^\varepsilon| \xrightarrow{P} 0\ \text{ as } n\to\infty.	
$$
Therefore, we can ascertain
$$
	{P}\bigg(\sup_{0\leq t \leq 1}|Y_{t}^{n,\varepsilon} - S_t^0| > \gamma\bigg)\to 0\text{ as } \varepsilon\to 0,\,n\to\infty.\ \ \ \ \ \ \ \ \ \ \ \ \ \ \ \ \ \ \ \ \ \ \square
$$

\begin{lem}\label{lemma2}
	Suppose {\bf (A1)} and {\bf(A2)} hold. Let $f \in C^{1,1,1}_\uparrow\left( [0,1] \times \mathbb{R}^d \times \Theta \right)$. Then,
$$
	\frac{1}{n}\sum\limits_{k=1}^n f(t_{k-1},S_{t_{k-1}}^\varepsilon,\theta) \xrightarrow{P_{\theta_0}} \dint_0^1 f(s,S_s^0,\theta)ds.
$$
as $\varepsilon \to 0$ and $n \to \infty$, uniformly in $\theta \in \Theta$.
\end{lem}
{\bf Proof.}\ \ We follow the argument of \cite[Lemma 3.3]{Sun}.  Using the condition imposed on $f$ and Lemma \ref{lem1}, we get
\begin{eqnarray*}
	&&\sup_{\theta \in \Theta} \bigg \lvert \frac{1}{n}\sum\limits_{k=1}^n f(t_{k-1},S_{t_{k-1}}^\varepsilon,\theta) - \nonumber
	\dint_0^1 f(s,S_s^0,\theta)ds \bigg \rvert \\ \nonumber
	& \leq&  \sup_{\theta \in \Theta} \sum_{k=1}^n \dint_{t_{k-1}}^{t_k}|f(t_{k-1},S_{t_{k-1}}^\varepsilon,\theta)  -   f(s,Y_s^{n,\varepsilon},\theta) |ds+   \sup_{\theta \in \Theta} \dint_0^1| f(s,Y_s^{n,\varepsilon},\theta) - f(s,S_s^0,\theta)| ds
 \\ \nonumber
 	& \leq&  \sum_{k=1}^n \dint_{t_{k-1}}^{t_k}\bigg (\int_0^1 \sup_{\theta \in \Theta} |(\nabla_t f) (t_{k-1} + v(s-t_{k-1}),Y_s^{n,\varepsilon},\theta)|dv \bigg )|s-t_{k-1}|ds\\ \nonumber
  &&+\dint_0^1 \bigg(\dint_0^1  \sup_{\theta \in \Theta} \nonumber
		\lvert \left( \nabla_x f \right) \left(s, S_s^0 + u \left( Y_s^{n,\varepsilon} -S_s^0 \right),\theta\right) \rvert du\bigg)\lvert Y_s^{n,\varepsilon} -S_s^0  \rvert  ds\\ \nonumber
	 	& \leq&  C \left( 1 + \sup_{s \in [0,1]} \lvert S_s^0 \rvert + \sup_{s \in [0,1]} \lvert S_s^{\varepsilon} \rvert \right)^{\lambda} \left[\frac{1}{n}+\sup_{s \in [0,1]}\lvert Y_s^{n,\varepsilon} -S_s^0  \rvert\right]\\
& \xrightarrow{P_{\theta_0}}&0 \text{ as } \varepsilon \to 0,\, n\to\infty.\ \ \ \ \ \ \ \ \ \ \ \ \ \ \ \ \ \ \ \ \ \ \ \ \ \ \ \ \ \ \ \ \ \ \ \ \ \ \ \ \ \ \ \ \ \ \ \ \ \ \ \ \square
\end{eqnarray*}

Define
$$
\tau_m^{\varepsilon} = \inf \{t \ge 0 : |S_t^\varepsilon| \vee |S_{t-}^\varepsilon| \ge m\}, \ \ \ \
\tau_m^0= \inf \{t \ge 0 :  |S_t^0| \ge m\}.
$$
Similar to \cite[Lemma 4.3]{Long}, we can prove the following lemma.
\begin{lem}\label{lem3}
	Suppose {\bf (A1)} and {\bf(A2)} hold. Then,   for any
	$m >0$, $\tau_m^{\varepsilon}  \xrightarrow{P_{\theta_0}} \tau_m^0$ as $\varepsilon\to 0$.
\end{lem}

\begin{lem}\label{lem5}
Suppose {\bf (A1)}--{\bf(A3)} hold. 	Let $f \in C^{1,1,1}_\uparrow\left( [0,1] \times \mathbb{R}^d \times \Theta \right)$. Then, for $1 \leq i \leq d$,
$$
	\sum\limits_{k=1}^n f(t_{k-1},S_{t_{k-1}}^\varepsilon,\theta) \left(S_{t_{k}}^{\varepsilon,i} - S_{t_{k-1}}^{\varepsilon,i} -b_i(t_{k-1},S_{t_{k-1}}^{\varepsilon}, \theta_0)\Delta_{t_{k-1}}  \right) \xrightarrow{P_{\theta_0}} 0
$$
as $\varepsilon\to 0$ and  $n\to\infty$, uniformly in $\theta \in \Theta$, where $S_t^{\varepsilon,i}$ and $b_i$ are the $i$-{th} components of $S^{\varepsilon}_t$ and $b$, respectively, and $\Delta_{t_{k-1}}:=t_k-t_{k-1}=\frac{1}{n}$.
\end{lem}
{\bf Proof.}\ \ Recall that
\begin{eqnarray*}
	S_{t_k}^{\varepsilon, i} &=& S_{t_{k-1}}^{\varepsilon, i} + \dint_{t_{k-1}}^{t_{k}} b_i(s,S_s^{\varepsilon},\theta_0)ds
 + \varepsilon\dint_{t_{k-1}}^{t_{k}} \bigg (\sigma_i(s,S_s^{\varepsilon} )dB_s  \\ && + \dint_{\{|u| \leq 1\}} H_i(s,S_{s-}^{\varepsilon},u)\widetilde{N}(ds,du)
  + \dint_{\{|u| > 1\}} G_i(s,S_{s-}^{\varepsilon},u)N(ds,du) \bigg ).
\end{eqnarray*}
Then, we have
\begin{eqnarray*}
	&&\sum\limits_{k=1}^n f(t_{k-1},S_{t_{k-1}}^\varepsilon,\theta) \left(S_{t_{k}}^{\varepsilon,i} - S_{t_{k-1}}^{\varepsilon,i} -b_i(t_{k-1},S_{t_{k-1}}^{\varepsilon}, \theta_0)\Delta_{t_{k-1}}  \right)\\
	& =& \sum\limits_{k=1}^n\dint_{t_{k-1}}^{t_{k}}  f(t_{k-1},S_{t_{k-1}}^{\varepsilon},\theta) \left( b_i(s,S_s^{\varepsilon}, \theta_0)-b_i(t_{k-1},S_{t_{k-1}}^{\varepsilon}, \theta_0) \right)ds \\ \nonumber
	&& + \varepsilon \sum_{k=1}^n\dint_{t_{k-1}}^{t_{k}} f(t_{k-1},S_{t_{k-1}}^{\varepsilon},\theta)\sigma_i(s,S_{s}^{\varepsilon} )dB_s\\ \nonumber
	&& + \varepsilon \sum_{k=1}^n\dint_{t_{k-1}}^{t_{k}} f(t_{k-1},S_{t_{k-1}}^{\varepsilon},\theta)\dint_{\{|u| \leq 1\}} H_i(s,S_{s-}^{\varepsilon},u)\widetilde{N}(ds,du) \\ \nonumber
	&& + \varepsilon \sum_{k=1}^n\dint_{t_{k-1}}^{t_{k}} f(t_{k-1},S_{t_{k-1}}^{\varepsilon},\theta) \dint_{\{|u| > 1\}} G_i(s,S_{s-}^{\varepsilon},u)N(ds,du)
\end{eqnarray*}
\begin{eqnarray*}
	& =&\dint_0^1  f(s,Y_s^{n,\varepsilon},\theta) \left( b_i(s,S_s^{\varepsilon}, \theta_0)-b_i(s,Y_s^{n,\varepsilon}, \theta_0) \right)ds \\ \nonumber
	&& + \varepsilon\dint_0^1  f(s,Y_s^{n,\varepsilon},\theta)\sigma_i(s,S_{s}^{\varepsilon} )dB_s\\ \nonumber
	&& + \varepsilon\dint\limits_0^1  f(s,Y_s^{n,\varepsilon},\theta)\dint_{\{|u| \leq 1\}}  H_i(s,S_{s-}^{\varepsilon},u)\widetilde{N}(ds,du) \\ \nonumber
	&& + \varepsilon\dint\limits_0^1  f(s,Y_s^{n,\varepsilon},\theta) \dint_{\{|u| > 1\}}  G_i(s,S_{s-}^{\varepsilon},u)N(ds,du)\\
&&+\sum\limits_{k=1}^n\dint_{t_{k-1}}^{t_{k}} [ f(t_{k-1},S_{t_{k-1}}^{\varepsilon},\theta) -f(s,Y_s^{n,\varepsilon},\theta)]\left(b_i(s,S_s^{\varepsilon}, \theta_0) -b_i(s,Y_s^{n,\varepsilon}, \theta_0) \right)ds\\
&&+\sum\limits_{k=1}^n\dint_{t_{k-1}}^{t_{k}}  f(t_{k-1},S_{t_{k-1}}^{\varepsilon},\theta) \left( b_i(s,Y_s^{n,\varepsilon}, \theta_0)-b_i(t_{k-1},Y_s^{n,\varepsilon}, \theta_0) \right)ds \\ \nonumber
	&& + \varepsilon \sum_{k=1}^n\dint_{t_{k-1}}^{t_{k}} [f(t_{k-1},S_{t_{k-1}}^{\varepsilon},\theta)- f(s,Y_s^{n,\varepsilon},\theta)]\sigma_i(s,S_{s}^{\varepsilon} )dB_s\\ \nonumber
	&& + \varepsilon \sum_{k=1}^n\dint_{t_{k-1}}^{t_{k}} [f(t_{k-1},S_{t_{k-1}}^{\varepsilon},\theta)- f(s,Y_s^{n,\varepsilon},\theta)]\dint_{\{|u| \leq 1\}} H_i(s,S_{s-}^{\varepsilon},u)\widetilde{N}(ds,du) \\ \nonumber
	&& + \varepsilon \sum_{k=1}^n\dint_{t_{k-1}}^{t_{k}} [f(t_{k-1},S_{t_{k-1}}^{\varepsilon},\theta) - f(s,Y_s^{n,\varepsilon},\theta)]\dint_{\{|u| > 1\}} G_i(s,S_{s-}^{\varepsilon},u)N(ds,du).
\end{eqnarray*}
Using the condition imposed on $f$, (A2), (A3), Lemma \ref{lem1},  Lemma \ref{lem3} and following the argument as in the proof of \cite[Lemma 3.5]{Sun}, we can show that all of the above terms converge to zero in probability as $\varepsilon\to 0$ and  $n\to\infty$, uniformly in $\theta \in \Theta$. \hfill $\square$

\vskip 0.5cm

\noindent {\bf Proof of Theorem \ref{thm1}.}\ \ Define
$$
\Phi_{n,\varepsilon}(\theta):=\Psi_{n,\varepsilon}(\theta)-\Psi_{n,\varepsilon}(\theta_0).
$$
We have
$$
\hat{\theta}_{n,\varepsilon}= \arg\min_{\theta\in\Theta}\Phi_{n,\varepsilon}(\theta),
  $$
and
\begin{eqnarray*}
	\Phi_{n,\varepsilon}(\theta) &=& -2 \sum\limits_{k=1}^n (b(t_{k-1},S_{t_{k-1}}^\varepsilon,\theta) - b(t_{k-1},S_{t_{k-1}}^\varepsilon,\theta_0))^\top\left(S_{t_k} ^\varepsilon-
		S_{t_{k-1}}^\varepsilon- \frac{1}{n}b(t_{k-1},S_{t_{k-1}}^\varepsilon,\theta_0)\right) \\
	&& + \frac{1}{n}\sum\limits_{k=1}^n \lvert b(t_{k-1},S_{t_{k-1}}^\varepsilon,\theta) - b(t_{k-1},S_{t_{k-1}}^\varepsilon,\theta_0) \rvert^2\\
&:=&	\Phi^{(1)}_{n,\varepsilon}(\theta)	+		\Phi^{(2)}_{n,\varepsilon}(\theta).
\end{eqnarray*}

By Lemma \ref{lem5} and letting
$f(t,x,\theta) = b_i(t,x,\theta)-b_i(t,x,\theta_0)$, $1 \leq i \leq d$,
we get $\sup_{\theta \in \Theta}|\Phi^{(1)}_{n,\varepsilon}(\theta)| \xrightarrow{{P}_{\theta_0}} 0$ as $\varepsilon \to 0$ and $n \to \infty$. By Lemma \ref{lemma2} and letting
$f(t,x,\theta) = |b(t,x,\theta)-b(t,x,\theta_0)|^2$, we get $\sup_{\theta \in \Theta}|\Phi^{(2)}_{n,\varepsilon}(\theta)-F(\theta)| \xrightarrow{{P}_{\theta_0}} 0$ as $\varepsilon \to 0$ and $n \to \infty$, where
$$
F(\theta) := \dint_0^1 |b(t,S_t^0,\theta)-b(t,S_t^0,\theta_0)|^2dt.
$$
Then, we have that
$$
	\sup_{\theta\in\Theta}\lvert \Phi_{n,\varepsilon}(\theta) - F(\theta)\rvert \xrightarrow{{P}_{\theta_0}} 0\ \ {\rm as}\ \varepsilon\to0,\,n\to\infty.
$$
By {\bf (A4)} and  the continuity of $S^0$, we get
$$
\sup_{|\theta - \theta_0|>\delta} -F(\theta) < -F(\theta_0) = 0,\ \ \ \  \forall \delta >0.
$$
Therefore, the proof is complete by \cite[Theorem 5.7]{van}. \hfill$\square$

\subsection{Proof of Theorem \ref{thm2}}

Note that
$$
	\nabla_\theta {\Phi}_{n,\varepsilon} = -2\sum_{k=1}^n (\nabla_\theta b)^\top(t_{k-1},S_{t_{k-1}}^\varepsilon,\theta)(S_{t_{k}}^\varepsilon-S_{t_{k-1}}^\varepsilon-b(t_{k-1},S_{t_{k-1}}^\varepsilon,\theta)\Delta_{t_{k-1}}).
$$
Define $G_{n,\varepsilon}(\theta) = (G_{n,\varepsilon}^i)_{1\leq i \leq p}^\top $, where
$$
	G_{n,\varepsilon}^i(\theta) := \sum_{k=1}^n(\partial_{\theta_i} b)^\top(t_{k-1},S_{t_{k-1}}^\varepsilon,\theta)(S_{t_{k}}^\varepsilon-S_{t_{k-1}}^\varepsilon-b(t_{k-1},S_{t_{k-1}}^\varepsilon,\theta)\Delta_{t_{k-1}}),\ \ \ \ 1\le i\le p,
$$
and define
$$
	K_{n,\varepsilon}(\theta) := \nabla_\theta G_{n,\varepsilon}(\theta),$$
namely, a $p\times p$ matrix comprised of $(\partial_{\theta_j}G^i_{n,\varepsilon})_{1\leq i,j \leq p}$. Lastly,
define
$$
	K^{ij}(\theta) =  \dint\limits_0^1 (\partial_{\theta_j}\partial_{\theta_i}b)^\top(s,S_s^0,\theta)(b(s,S_s^0,\theta_0)
	-b(s,S_s^0,\theta))ds- I^{ij} (\theta),
$$
and
$$
	K(\theta) = (K^{ij}(\theta))_{1\leq i,j \leq p}.
$$

\begin{lem}\label{lem4}
Suppose {\bf (A1)}, {\bf(A2)} and {\bf(A5)} hold. 	Let $f \in C^{1,1,1}_\uparrow\left( [0,1] \times \mathbb{R}^d \times \Theta \right)$. Then, for $1 \leq i \leq d$ and each $\theta \in \Theta$,
\begin{eqnarray*}
&&\sum_{k=1}^n f(t_{k-1},S_{t_{k-1}}^{\varepsilon},\theta)\bigg( \int_{t_{k-1}}^{t_{k}} \sigma(s,S_s^\varepsilon)dB_s
+\int_{t_{k-1}}^{t_{k}}\int_{\{|u|\le1\}}H_i(s,S_{s-}^\varepsilon,u)\widetilde{N}
(ds,du)\\ \nonumber
&&\ \ \ \ \ +\int_{t_{k-1}}^{t_{k}}\int_{\{|u|> 1\}}G_i(s,S_{s-}^\varepsilon,u)N
(ds,du)\bigg)  \\ \nonumber
&\xrightarrow{P_{\theta_0}}&\dint_0^1  f(s,S_s^0,\theta)\bigg(  \sigma_i(s,S_s^0)dB_s+\int_{\{|u|\le1\}}H_i(s,S_{s}^0,u)\widetilde{N}
(ds,du)\\ \nonumber
&&\ \ \ \ \ +\int_{\{|u|> 1\}}G_i(s,S_{s}^0,u)N
(ds,du)\bigg),
\end{eqnarray*}
as $\varepsilon \to 0$ and $n\to \infty$.
\end{lem}

\noindent {\bf Proof.}\ \ We follow the arguments of \cite[Lemma 3.4]{Sun}  and \cite[Lemma 4.6]{Long}. Define
\begin{eqnarray*}
A_{n,\varepsilon}&:=&\sum_{k=1}^n f(t_{k-1},S_{t_{k-1}}^{\varepsilon},\theta)\bigg( \int_{t_{k-1}}^{t_{k}} \sigma_i(s,S_{s}^\varepsilon)dB_s
+\int_{t_{k-1}}^{t_{k}}\int_{\{|u|\le1\}}H_i(s,S_{s-}^\varepsilon,u)\widetilde{N}
(ds,du)\\
&&\ \ \ \ \ +\int_{t_{k-1}}^{t_{k}}\int_{\{|u|> 1\}}G_i(s,S_{s-}^\varepsilon,u)N
(ds,du)\bigg), \\
B_{n,\varepsilon}&:=&\dint_0^1  f(s,Y_s^{n,\varepsilon},\theta)\bigg(  \sigma_i(s,S_s^\varepsilon)dB_s+\int_{\{|u|\le1\}}H_i(s,S_{s-}^\varepsilon,u)\widetilde{N}
(ds,du)\\
&&\ \ \ \ \ +\int_{\{|u|> 1\}}G_i(s,S_{s-}^\varepsilon,u)N
(ds,du)\bigg).
\end{eqnarray*}
We can show that  $A_{n,\varepsilon}=B_{n,\varepsilon}+o_{P}(1)$. In fact, we have
\begin{eqnarray*}
&&|A_{n,\varepsilon}-B_{n,\varepsilon}|\\
&\leq& \sum_{k=1}^n \int_{t_{k-1}}^{t_k} |f(t_{k-1},S_{t_{k-1}}^\varepsilon,\theta)-f(s,Y_s^{n,\varepsilon},\theta)|\int_{\{|u|> 1\}}|G_i(s,S_{s-}^\varepsilon,u)|N
(ds,du) \\
&&+  \bigg| \sum_{k=1}^n \int_{t_{k-1}}^{t_k} [f(t_{k-1},S_{t_{k-1}}^\varepsilon,\theta)-f(s,Y_s^{n,\varepsilon},\theta)]\bigg( \sigma_i(s,S_{s}^\varepsilon)dB_s+\int_{\{|u|\le1\}}H_i(s,S_{s-}^\varepsilon,u)\widetilde{N}
(ds,du)\bigg)\bigg|\\
&:=&I_{n,\varepsilon}+J_{n,\varepsilon},
\end{eqnarray*}
where
\begin{eqnarray*}
I_{n,\varepsilon}&\leq & \sup_{\theta \in \Theta} \sum_{k=1}^n \int_{t_{k-1}}^{t_k} \bigg(\int_0^1 |(\nabla_t f)(t_{k-1}+v(s-t_{k-1}),Y_s^{n,\varepsilon},\theta)|dv\bigg)|s-t_{k-1}|  \\\nonumber
&&\cdot\int_{\{|u|> 1\}}|G_{i}(s,S_{s-}^\varepsilon,u)|N(ds,du)\\\nonumber
&\leq &  \frac{CK}{n}\left(1+\sup_{t\in[0,1]} |S_t^{\varepsilon}|\right)^{\lambda+1} \int_{\{|u|> 1\}}\xi(u)N(ds,du)\\
& \xrightarrow{P}& 0\ {\rm as}\ \varepsilon \to0,\,n\to\infty.
\end{eqnarray*}
For $J_{n,\varepsilon}$ and $\upeta> 0$, using the stopping time $\tau_m^\varepsilon$, Lemma \ref{lem1}, Markov's inequality and dominated convergence, we get
\begin{eqnarray*}
&&{P}\bigg(\bigg| \sum_{k=1}^n\int_{t_{k-1}}^{t_k} [f(t_{k-1},S_{t_{k-1}}^\varepsilon,\theta)-f(s,Y_s^{n,\varepsilon},\theta)]\bigg( \sigma_i(s,S_{s}^\varepsilon){1}_{\{s\leq\tau^{\varepsilon}_m \}}dB_s\\
&&\ \ \ \ \ \ \ \ \ \ \ \ \ \ \ \ \ +\int_{\{|u|\le1\}}H_i(s,S_{s-}^\varepsilon,u){1}_{\{s\leq\tau^{\varepsilon}_m \}}\widetilde{N}
(ds,du)\bigg)\bigg| > \upeta\bigg)\\
&\leq &\frac{1}{\upeta} \bigg({E}\left\{\sum_{k=1}^n\int_{t_{k-1}}^{t_k}|f(t_{k-1},S_{t_{k-1}}^\varepsilon,\theta)-f(s,Y_s^{n,\varepsilon},\theta)|^2\sigma_i^2(s,S_{s}^\varepsilon)
{1}_{\{s\leq\tau^{\varepsilon}_m \}}ds\right\}\bigg)^{1/2} \\\nonumber
&&
+ \frac{1}{\upeta} \bigg({E}\left\{\sum_{k=1}^n\int_{t_{k-1}}^{t_k}|f(t_{k-1},S_{t_{k-1}}^\varepsilon,\theta)-f(s,Y_s^{n,\varepsilon},\theta)|^2{1}_{\{s\leq\tau^{\varepsilon}_m\}}ds\right\}\bigg)^{1/2}\\
&&\ \ \ \ \ \cdot \bigg({E}\left\{\int_0^1			\dint_{\{|u| \le 1\}} |H_i(s,S_{s-}^\varepsilon,u)|^2 \mu(du){1}_{\{s\leq\tau^{\varepsilon}_m\}}\right\}ds \bigg)^{1/2}\\
&\leq &\frac{1}{\upeta} \bigg({E}\left\{\sum_{k=1}^n\int_{t_{k-1}}^{t_k}\bigg(\int_0^1 |(\nabla_t f)(t_{k-1}+v(s-t_{k-1}),Y_s^{n,\varepsilon},\theta)|dv|s-t_{k-1}|\bigg)^2\sigma_i^2(s,S_{s}^\varepsilon)
{1}_{\{s\leq\tau^{\varepsilon}_m \}}ds\right\}\bigg)^{1/2} \\\nonumber
&&
+ \frac{1}{\upeta} \bigg({E}\left\{ \sum_{k=1}^n\int_{t_{k-1}}^{t_k}\bigg(\int_0^1 |(\nabla_t f)(t_{k-1}+v(s-t_{k-1}),Y_s^{n,\varepsilon},\theta)|dv|s-t_{k-1}|\bigg)^2{1}_{\{s\leq\tau^{\varepsilon}_m\}} ds\right\}\bigg)^{1/2}\\
&&\ \ \ \ \ \cdot \bigg({E}\left\{\int_0^1			\dint_{\{|u| \le 1\}} |H_i(s,S_{s-}^\varepsilon,u)|^2 \mu(du){1}_{\{s\leq\tau^{\varepsilon}_m\}}ds \right\}\bigg)^{1/2}\\
&\leq &\frac{1}{n\upeta}\bigg({E}\left\{\sum_{k=1}^n\int_{t_{k-1}}^{t_k}\bigg(C(1+|Y_s^{n,\varepsilon}|)^\lambda\bigg)^2\sigma_i^2(s,S_{s}^\varepsilon)
{1}_{\{s\leq\tau^{\varepsilon}_m \}}ds\right\}\bigg)^{1/2} \\\nonumber
&&
+ \frac{K}{n\upeta}\bigg({E}\left\{ \sum_{k=1}^n\int_{t_{k-1}}^{t_k}\bigg(C(1+|Y_s^{n,\varepsilon}|)^\lambda\bigg)^2{1}_{\{s\leq\tau^{\varepsilon}_m\}}ds \right\}\bigg)^{1/2}\\
&&\ \ \ \ \ \ \  \cdot \bigg({E}\left\{\int_0^1(1+|Y_s^{n,\varepsilon}|)^2{1}_{\{s\leq\tau^{\varepsilon}_m\}}ds\right\}\bigg)^{1/2}\bigg(\int_{\{|u|\le 1\}}\eta^2(u)\mu(du)\bigg)^{1/2}\\
&\rightarrow&0\ {\rm as}\ \varepsilon \to0,\,n\to\infty,
\end{eqnarray*}
and
\begin{eqnarray*}
&&{P}\bigg(  \bigg|\sum_{k=1}^n\int_{t_{k-1}}^{t_k}[f(t_{k-1},S_{t_{k-1}}^\varepsilon,\theta)-f(s,Y_s^{n,\varepsilon},\theta)]\bigg( \sigma_i(s,S_{s}^\varepsilon)dB_s\\
&&\ \ \ \ \ \ \ \ \ \ \ +\int_{\{|u|\le1\}}H_i(s,S_{s-}^\varepsilon,u)\widetilde{N}
(ds,du)\bigg)\bigg| > \upeta\bigg) \\\nonumber
&\leq &{P}\bigg(\bigg| \sum_{k=1}^n\int_{t_{k-1}}^{t_k} [f(t_{k-1},S_{t_{k-1}}^\varepsilon,\theta)-f(s,Y_s^{n,\varepsilon},\theta)]\bigg( \sigma_i(s,S_{s}^\varepsilon){1}_{\{s\leq\tau^{\varepsilon}_m \}}dB_s\\
&&\ \ \ \ \ \ \ \ \ \ \  +\int_{\{|u|\le1\}}H_i(s,S_{s-}^\varepsilon,u){1}_{\{s\leq\tau^{\varepsilon}_m \}}\widetilde{N}
(ds,du)\bigg)\bigg| > \upeta\bigg)  + {P}(\tau_m^{\varepsilon}<1)\\
&\rightarrow&0\ {\rm as}\ \varepsilon \to0,\,n\to\infty,
\end{eqnarray*}
by Lemma \ref{lem3}.

Finally, using the condition imposed on $f$, (A5), \cite[Problem 13, page 151]{ei}, the continuous mapping theorem, Lemma \ref{lem1} and similar arguments as above, we can show that
\begin{eqnarray*}
B_{n,\varepsilon}
&\xrightarrow{P_{\theta_0}}&\dint_0^1  f(s,S_s^0,\theta)\bigg(  \sigma_i(s,S_s^0)dB_s+\int_{\{|u|\le1\}}H_i(s,S_{s}^0,u)\widetilde{N}
(ds,du)\\ \nonumber
&&\ \ \ \ \ +\int_{\{|u|> 1\}}G_i(s,S_{s}^0,u)N
(ds,du)\bigg)
\end{eqnarray*}
as $\varepsilon \to 0$ and $n\to \infty$.\hfill $\square$

\begin{lem}\label{lem7}
	Suppose {\bf (A1)}--{\bf(A6)} hold. Then, for $1\leq i \leq p$,
\begin{eqnarray*}
		\varepsilon^{-1}G_{n,\varepsilon}^i(\theta_0) &\xrightarrow{P_{\theta_0}}& \int_0^1 (\partial_{\theta_i}b)^\top(s,S_s^0,\theta_0)\bigg(  \sigma(s,S_s^0)dB_s+\int_{\{|u|\le1\}}H(s,S_{s}^0,u)\widetilde{N}
(ds,du)\\
&&\ \ \ \ \ +\int_{\{|u|> 1\}}G(s,S_{s}^0,u)N
(ds,du)\bigg)
\end{eqnarray*}
as $\varepsilon \to 0$ and $n\to \infty$.
\end{lem}
\noindent {\bf Proof.}\ \ We follow the arguments of \cite[Lemma 3.6]{Sun}  and \cite[Lemma 4.7]{Long}.  For $1 \leq i \leq p$, we have
\begin{eqnarray*}
&&\varepsilon^{-1}G_{n,\varepsilon}^i(\theta_0)\\
 &=& \varepsilon^{-1}\sum_{k=1}^n(\partial_{\theta_i}b)^\top({t_{k-1}},S_{t_{k-1}}^\varepsilon,\theta_0)(S_{t_{k}}^\varepsilon-S_{t_{k-1}}^\varepsilon-b(t_{k-1},S_{t_{k}}^\varepsilon,\theta_0)\Delta_{k-1})\\ \nonumber
 	&=& \varepsilon^{-1}\sum_{k=1}^n(\partial_{\theta_i}b)^\top({t_{k-1}},S_{t_{k-1}}^\varepsilon,\theta_0) \int_{t_{k-1}}^{t_{k}}(b(s,S_{s}^\varepsilon,\theta_0)-b(t_{k-1},S_{t_{k}}^\varepsilon,\theta_0))ds \\ \nonumber
&&+ \sum_{k=1}^n(\partial_{\theta_i}b)^\top({t_{k-1}},S_{t_{k-1}}^\varepsilon,\theta_0)\int_{t_{k-1}}^{t_{k}}\bigg(  \sigma(s,S_s^\varepsilon)dB_s+\int_{\{|u|\le1\}}H(s,S_{s-}^\varepsilon,u)\widetilde{N}(ds,du)\\ \nonumber
&&\ \ \ \ \ \ \ \ +\int_{\{|u|> 1\}}G(s,X_{s-}^\varepsilon,u)N
(ds,du)\bigg) \\ \nonumber
&:=& H^{(1)}_{n,\varepsilon}(\theta_0)+H^{(2)}_{n,\varepsilon}(\theta_0).
\end{eqnarray*}
Using Lemma \ref{lem4} with $f(t,x,\theta)=(\partial_{\theta_i}b_j(t,x,\theta))^\top(t,x,\theta)$
and $\theta=\theta_0$ for $1 \leq i \leq p$, $1\leq j \leq d $, we get
\begin{eqnarray*}
	H^{(2)}_{n,\varepsilon}(\theta_0)& \xrightarrow{{P}_{\theta_0}}& \int_0^1 (\partial_{\theta_i}b)^\top(s,S_s^0,\theta_0)\bigg(  \sigma(s,S_s^0)dB_s+\int_{\{|u|\le1\}}H(s,S_{s}^0,u)\widetilde{N}
(ds,du)\\ \nonumber
&&+\int_{\{|u|> 1\}}G(s,S_{s}^0,u)N
(ds,du)\bigg)
\end{eqnarray*}
as $\varepsilon\to0$ and $n\to\infty$.

For $H^{(1)}_{n,\varepsilon}(\theta_0)$, given $s \in [t_{k-1},t_k]$, we have that
\begin{eqnarray*}
	&&S_s^\varepsilon-S_{t_{k-1}}^\varepsilon\\
&= &\dint_{t_{k-1}}^s(b(r,S_r^\varepsilon,\theta_0)-b(t_{k-1},S_{t_{k-1}}^\varepsilon,\theta_0))dr + b(t_{k-1},S_{t_{k-1}}^\varepsilon,\theta_0)(s-t_{k-1}) \\ \nonumber
&&+
	\varepsilon \dint_{t_{k-1}}^s\bigg(  \sigma(s,S_s^\varepsilon)dB_s+\int_{\{|u|\le1\}}H(s,S_{s-}^\varepsilon,u)\widetilde{N}(ds,du)+\int_{\{|u|> 1\}}G(s,S_{s-}^\varepsilon,u)N
(ds,du)\bigg).
\end{eqnarray*}
Using the Lipschitz condition on $b$ and the Cauchy-Schwarz inequality, we get
\begin{eqnarray*}
	|S_s^\varepsilon-S_{t_{k-1}}^\varepsilon|^2 & \leq& 2\bigg| \dint_{t_{k-1}}^s(b(r,S_r^\varepsilon,\theta_0)-b(t_{k-1},S_{t_{k-1}}^\varepsilon,\theta_0))dr \bigg|^2 \\ \nonumber &&+ \nonumber
	2\bigg\{|b(t_{k-1},S_{t_{k-1}}^\varepsilon,\theta_0)|(s-t_{k-1}) +\varepsilon \dint_{t_{k-1}}^s\bigg(  \sigma(r,S_r^\varepsilon)dB_r \\ \nonumber
&&+\int_{\{|u|\le1\}}H(r,S_{r-}^\varepsilon,u)\widetilde{N}(dr,du)+\int_{\{|u|> 1\}}G(r,S_{r-}^\varepsilon,u)N(dr,du)\bigg) \nonumber \bigg\}^2\\
	& \leq& 2K^2n^{-1}\bigg(\dint_{t_{k-1}}^s  |S_r^\varepsilon-S_{t_{k-1}}^\varepsilon|dr\bigg)^2 + 2\bigg\{n^{-1}|b(t_{k-1},S_{t_{k-1}}^\varepsilon,\theta_0)| \\ \nonumber
	&& + \varepsilon  \sup_{{t_{k-1}} \leq s \leq t_k} \bigg|\dint_{t_{k-1}}^s\bigg(  \sigma(r,S_r^\varepsilon)dB_r +\int_{\{|u|\le1\}}H(r,S_{r-}^\varepsilon,u)\widetilde{N}(dr,du)
	\\ \nonumber
&&+\int_{\{|u|> 1\}}G(r,S_{r-}^\varepsilon,u)N(dr,du)\bigg) \bigg| \bigg\}^2.
\end{eqnarray*}
Now applying Gronwall's inequality, we get
\begin{eqnarray*}
	&&\sup_{t_{k-1}\leq s \leq t_k} |S_s^\varepsilon-S_{t_{k-1}}^\varepsilon|\\
   &\leq& \sqrt{2}e^{K^2n^{-2}}\bigg(n^{-1}|b(t_{k-1},S_{t_{k-1}}^\varepsilon,\theta_0)| \\ \nonumber
	&&+ \varepsilon\sup_{t_{k-1}\leq s \leq t_k}\bigg\vert\dint_{t_{k-1}}^s \sigma(r,S_r^\varepsilon)dB_r +\int_{\{|u|\le1\}}H(r,S_{r-}^\varepsilon,u)\widetilde{N}(dr,du)	+\int_{\{|u|> 1\}}G(r,S_{r-}^\varepsilon,u)N(dr,du)\bigg\vert \bigg).
\end{eqnarray*}
Further, by the Lipschitz condition on $b$ and (A3), we obtain that
\begin{eqnarray*}
	&&|H^{(1)}_{n,\varepsilon}(\theta_0)|{1}_{\{1\leq\tau^{\varepsilon}_m \}}\\
	&\leq&
	\varepsilon^{-1}\sum_{k=1}^n|(\partial_{\theta_i}b)^\top({t_{k-1}},S_{t_{k-1}}^\varepsilon,\theta_0)|\cdot \nonumber
	\int_{t_{k-1}}^{t_{k}}|(b(s,S_{s}^\varepsilon,\theta_0)-b(t_{k-1},S_{t_{k}}^\varepsilon,\theta_0))|{1}_{\{s\leq\tau^{\varepsilon}_m \}}ds \\
&	\leq&
	\varepsilon^{-1}\sum_{k=1}^n|(\partial_{\theta_i}b)^\top({t_{k-1}},S_{t_{k-1}}^\varepsilon,\theta_0)|\cdot \nonumber
	\left[\int_{t_{k-1}}^{t_{k}}K|S_s^\varepsilon - S_{t_{k-1}}^\varepsilon|ds+Cn^{-1}(1+m)^{\lambda}\right]  \\
	&\leq& (n\varepsilon)^{-1}\sum_{k=1}^n|(\partial_{\theta_i}b)^\top({t_{k-1}},S_{t_{k-1}}^\varepsilon,\theta_0)|\left[K|\sup_{t_{k-1}\leq s \leq t_k}|S_s^\varepsilon - S_{t_{k-1}}^\varepsilon| +C(1+m)^{\lambda}\right]
\end{eqnarray*}
\begin{eqnarray*}
	&\leq& \sqrt{2}(n\varepsilon)^{-1}Ke^{K^2n^{-2}}n^{-1}\sum_{k=1}^n|(\partial_{\theta_i}b)^\top({t_{k-1}},S_{t_{k-1}}^\varepsilon,\theta_0)|\cdot
	|b({t_{k-1}},S_{t_{k-1}}^\varepsilon,\theta_0)| \\ \nonumber
&& +\sqrt{2}
	n^{-1}Ke^{K^2n^{-2}}\sum_{k=1}^n|(\partial_{\theta_i}b)^\top({t_{k-1}},S_{t_{k-1}}^\varepsilon,\theta_0)| \\ \nonumber
&&\ \ \cdot\sup_{t_{k-1}\leq s \leq t_k}\bigg\vert\dint_{t_{k-1}}^s \sigma(r,S_r^\varepsilon)dB_r
	+\int_{\{|u|\le1\}}H(r,S_{r-}^\varepsilon,u)\widetilde{N}(dr,du)
	+\int_{\{|u|> 1\}}G(r,S_{r-}^\varepsilon,u)N(dr,du)\bigg\vert \\
&&+ C(1+m)^{\lambda}(n\varepsilon)^{-1}\sum_{k=1}^n|(\partial_{\theta_i}b)^\top({t_{k-1}},S_{t_{k-1}}^\varepsilon,\theta_0)|\\
&:=&H^{(1,1)}_{n,\varepsilon}(\theta_0) + H^{(1,2)}_{n,\varepsilon}(\theta_0)+ H^{(1,3)}_{n,\varepsilon}(\theta_0).
\end{eqnarray*}
By (A2) and (A3), we can show that all of the above terms converge to zero in probability as $\varepsilon\to 0$ and  $n\to\infty$. Finally, the proof is complete by Lemma \ref{lem3}. \hfill $\square$

\begin{lem}\label{lem8}
	Assume that conditions (A1)-(A6) hold and $I(\theta_0)$ is positive definite. Then,
$$
		\sup_{\theta\in\Theta}|K_{n,\varepsilon}(\theta) - K(\theta)| \xrightarrow{{P}_{\theta_0}} 0
$$
	as $\varepsilon\to0$ and $n\to\infty$.
\end{lem}
{\bf Proof.}\ \ We follow the argument of \cite[Lemma 3.7]{Sun}.
For $1\leq i,j \leq p$, we have
\begin{eqnarray*}
	K_{n,\varepsilon}^{ij}(\theta)&=&\partial_{\theta_j}G^i_{n,\varepsilon}(\theta)\\
	&=& \sum_{k=1}^n(\partial_{\theta_j}\partial_{\theta_i}b)^\top(t_{k-1},S^\varepsilon_{t_{k-1}},\theta)
	(S_{t_{k}}^\varepsilon-S_{t_{k-1}}^\varepsilon-b(t_{k-1},S_{t_{k-1}}^\varepsilon,\theta_0)\Delta_{t_{k-1}})\\
	&&+n^{-1}\sum_{k=1}^n\bigg( (\partial_{\theta_j}\partial_{\theta_i}b)^\top(t_{k-1},S^\varepsilon_{t_{k-1}},\theta)
	(b(t_{k-1},S_{t_{k-1}}^\varepsilon,\theta_0)-b(t_{k-1},S_{t_{k-1}}^\varepsilon,\theta))\\
	&&\ \ \ \ \ \ \ \ \ \ \ \ \ \ \  - (\partial_{\theta_i}b)^\top(t_{k-1},S^\varepsilon_{t_{k-1}},\theta)
	\partial_{\theta_j}b(t_{k-1},S^\varepsilon_{t_{k-1}},\theta)\bigg)\\
&:=& K_{n,\varepsilon}^{(1)}(\theta)+K_{n,\varepsilon}^{(2)}(\theta).
\end{eqnarray*}
By application of Lemma \ref{lem5}  and setting $f(t,x,\theta) = \partial _{\theta_j}\partial_{\theta_i} b_l(t,x,\theta)$, $1 \leq i, j \leq p, 1 \leq l \leq d$, we get
$$
	\sup_{\theta\in\Theta}|K_{n,\varepsilon}^{ij,(1)}(\theta)| \xrightarrow{{P}_{\theta_0}} 0 \text{ as } \varepsilon\to0,n\to\infty.
$$
By application of Lemma \ref{lemma2} and setting
 $$f(t,x,\theta) = (\partial_{\theta_j}\partial_{\theta_i} b)^\top (t,x,\theta)(b(t,x, \theta_0) - b(t,x, \theta))-(\partial_{\theta_i}b)^\top (t,x,\theta)\partial_{\theta_j}b(t,x,\theta),$$ we get
$$
	\sup_{\theta\in\Theta}|K_{n,\varepsilon}^{ij,(2)}(\theta) - K^{ij}(\theta)| \xrightarrow{{P}_{\theta_0}} 0 \text{ as } \varepsilon\to0,n\to\infty.
$$
Therefore, the proof is complete.\hfill $\square$

\vskip 0.5cm
\noindent {\bf Proof of Theorem \ref{thm2}.}\ \ The proof is very similar to that given in \cite{Sun}. For the sake of completeness, we still include it here. We follow the arguments of \cite{Uch} and \cite{Sun}.  Consider the closed ball
$$
	B(\theta_0;\rho) := \{\theta : |\theta-\theta_0|\leq \rho\},\ \ \ \ \rho > 0.
$$
Using the established consistency of estimator $\hat{\theta}_{n,\varepsilon}$ from Theorem \ref{thm1}, we can find a sequence $\eta_{n,\varepsilon} \to 0$ as $\varepsilon\to0$ and $n\to\infty$
such that
$$
	B(\theta_0;\eta_{n,\varepsilon}) \subset \Theta \text{ and } {P}_{\theta_0}(\hat{\theta}_{n,\varepsilon} \in B(\theta_0;\eta_{n,\varepsilon})) \to 1.
$$

By application of Taylor's theorem when $\hat{\theta}_{n,\varepsilon} \in B(\theta_0;\eta_{n,\varepsilon})$, we get
$$
	D_{n,\varepsilon}T_{n,\varepsilon} = \varepsilon^{-1}(G_{n,\varepsilon}(\hat{\theta}_{n,\varepsilon})-G_{n,\varepsilon}(\theta_0)),
$$
where $D_{n,\varepsilon} := \int_0^1 K_{n,\varepsilon}(\theta_0 + u(\hat{\theta}_{n,\varepsilon} - \theta_0))du$ and
$T_{n,\varepsilon} := \varepsilon^{-1}(\hat{\theta}_{n,\varepsilon} - \theta_0) $ as $B(\theta_0;\eta_{n,\varepsilon})$ is a convex
subset of the parameter space $\Theta$. We have
\begin{eqnarray*}
	&&|D_{n,\varepsilon}-K_{n,\varepsilon}(\theta_0)|{1}_{\{\hat{\theta}_{n,\varepsilon} \in B(\theta_0;\eta_{n,\varepsilon})\}}\\
& \leq&
	\sup_{\theta \in B(\theta_0;\eta_{n,\varepsilon})} |K_{n,\varepsilon}(\theta)-K_{n,\varepsilon}(\theta_0)|\\
	&\leq& \sup_{\theta \in B(\theta_0;\eta_{n,\varepsilon})} |K_{n,\varepsilon}(\theta)-K(\theta)|
	+  \sup_{\theta \in B(\theta_0;\eta_{n,\varepsilon})}|K(\theta)-K(\theta_0)| + |K_{n,\varepsilon}(\theta_0)-K(\theta_0)|.
\end{eqnarray*}
Hence, by Lemma \ref{lem8}, we get
$$
	D_{n,\varepsilon} \to K_{n,\varepsilon}\xrightarrow{{P}_{\theta_0}} 0 \text{ as } \varepsilon\to0,n\to\infty.
$$

Note that $K(\theta)$ is continuous with respect to $\theta$ and $-K(\theta_0)= I(\theta_0)$ is positive definite; hence, there exists $\delta > 0$ such that
\begin{equation}\label{inf_delta}
	\inf_{|w| = 1} |K(\theta_0)w| > 2\delta.
\end{equation}
Choosing such a $\delta >0$, there exist $\varepsilon(\delta)>0$ and $N(\delta)>0$ such that for all $\varepsilon \in (0,\varepsilon(\delta)),\,n >N(\delta) $
we have
$$
	B(\theta_0;\eta_{n,\varepsilon}) \subset \Theta \text{ and } |K(\theta)-K(\theta_0)|<\frac{\delta}{2}\ {\rm for}\ \theta \in B(\theta_0;\eta_{n,\varepsilon}).
$$
Choose a $\delta > 0$ satisfying (\ref{inf_delta}) and set
$$
	\Gamma_{n,\varepsilon} = \bigg\{ \sup_{|\theta-\theta_0|\leq \eta_{n,\varepsilon}} |K_{n,\varepsilon}(\theta)-K(\theta_0)| < \frac{\delta}{2};
	\hat{\theta}_{n,\varepsilon} \in B(\theta_0;\,\eta_{n,\varepsilon}) \bigg\}.
$$
For $\varepsilon \in (0,\varepsilon(\delta))$ and $n > N(\delta)$, on the set $\Gamma_{n,\varepsilon}$ we have
\begin{eqnarray*}
	\sup_{|w|=1}|(D_{n,\varepsilon}-K(\theta_0))w|
& \leq& \sup_{|w|=1} \bigg|\bigg(D_{n,\varepsilon} -
	\int_0^1 K_{n,\varepsilon}(\theta_0 + u(\hat{\theta}_{n,\varepsilon} - \theta_0))du \bigg)w \bigg|\\
&&
	+ \sup_{|w|=1} \bigg|\bigg(
	\int_0^1 K_{n,\varepsilon}(\theta_0 + u(\hat{\theta}_{n,\varepsilon} - \theta_0))du - K(\theta_0)\bigg)w \bigg|\\
& \leq&
	\sup_{|\theta-\theta_0| \leq \eta_{n,\varepsilon}} |K_{n,\varepsilon}(\theta) - K(\theta)| + \frac{\delta}{2}\\
& <& \delta.
\end{eqnarray*}
Then, by (\ref{inf_delta}), on the set $\Gamma_{n,\varepsilon}$ we have that
$$
	\inf_{|w|=1}|D_{n,\varepsilon}w| \ge \inf_{|w|=1}|K(\theta_0)w| - \inf_{|w|=1}|(D_{n,\varepsilon}-K(\theta_0))w| > \delta > 0.
$$
Set
$\mathfrak{D}_{n,\varepsilon} = \{D_{n,\varepsilon} \text{ is invertible},\hat{\theta}_{n,\varepsilon} \in B(\theta_0;\eta_{n,\varepsilon})\}$. Thus, by Lemma \ref{lem8}, we get
$$
	{P}_{\theta_0}(\mathfrak{D}_{n,\varepsilon})\ge {P}_{\theta_0}(\Gamma_{n,\varepsilon})\to 1 \text{ as } \varepsilon\to 0,n\to\infty.
$$


Set
$$
	U_{n,\varepsilon} = D_{n,\varepsilon}{1}_{\mathfrak{D}_{n,\varepsilon}} + I_{p\times p}{1}_{\mathfrak{D}_{n,\varepsilon}^c},
$$
where $I_{p\times p}$ is the identity matrix.
We have
$$
	|U_{n,\varepsilon}-K(\theta_0)| \leq
	|D_{n,\varepsilon}-
	K(\theta_0)|{1}_{\mathfrak{D}_{n,\varepsilon}} + |I_{p\times p}-K(\theta_0)|{1}_{\mathfrak{D}_{n,\varepsilon}^c}
	\xrightarrow{{P}_{\theta_0}} 0
$$
since  ${P}_{\theta_0}(\mathfrak{D}_{n,\varepsilon}) \to 1$. Then, by Lemma \ref{lem7}, we obtain that
\begin{eqnarray*}
	T_{n,\varepsilon} &=& U_{n,\varepsilon}^{-1}D_{n,\varepsilon}T_{n,\varepsilon}{1}_{\mathfrak{D}_{n,\varepsilon}} +
	T_{n,\varepsilon}{1}_{\mathfrak{D}_{n,\varepsilon}^c}\\
&	= &U_{n,\varepsilon}^{-1}(-\varepsilon^{-1}(G_{n,\varepsilon}(\theta_0))){1}_{\mathfrak{D}_{n,\varepsilon}} +
	T_{n,\varepsilon}{1}_{\mathfrak{D}_{n,\varepsilon}^c}\\
	& \xrightarrow{P_{\theta_0}}&I^{-1}(\theta_0)\left(\int_0^1(\partial_{\theta_i}b)^\top(r,S^0_r,\theta)\Bigg\{\sigma(r,S^0_r)dB_r\right.\\
	&&\left.
	\ \ \ \ +{\int_{\{|u|\le1\}}}H(r,S^0_r,u)\widetilde{N}
	(dr,du)+{\int_{\{|u|> 1\}}}G(r,S^0_r,u)N
	(dr,du)\Bigg\}\right)^\top_{1\le i\le p}
\end{eqnarray*}
as $\varepsilon \to 0$ and $n \to \infty$. The proof is complete. \hfill $\square$

\section{Conclusion}\setcounter{equation}{0}

Our aim was to study the ability of estimating parameters governing a periodic transmission for a stochastic SIR model.
The model comes in one of two forms, population numbers and population proportions. Simulation studies were given for both forms and we are able to effectively estimate the period of a periodic transmission function.

The theoretical results from Section 2 and the proofs in Section 5 show that these estimation efforts hold in general. Namely, the results are suitable for other SDE applications beyond studying stochastic SIR models. The provided results for the estimators adds to the existing literature on SDEs driven by L\'{e}vy noises with small coefficient $\varepsilon$.

Algebraic manipulation to obtain LSEs was not a feasbile task; hence, the LS-GD algorithm was a useful tool for our work.
However, this algorithm is not optimal, and a future inquiry is to consider to what extent it may be optimized.

The work here lays a foundation for further work on parameter estimation of the USSIR model and SDEs driven by small L\'{e}vy noises. We will be interested in applying results of this paper
to real world data (e.g., COVID-19) with the inclusion of filtering to alleviate measurement errors.

\begin{center}
{\bf \footnotesize Acknowledgements}
\end{center}
{\scriptsize This work was partially supported by the Natural Sciences and Engineering Research Council of
Canada (No. 4394-2018). We wish to thank Calcul Qu\'{e}bec a regional partner of Digital Research Alliance of Canada for providing
high-performance computing
resources to accomplish the work presented here.}

\end{document}